\documentclass[12pt,a4paper]{article}

\usepackage{amssymb, amsmath, amsthm}

\usepackage[left=2.4cm,right=2.4cm,top=2.5cm,bottom=2.5cm,headheight=50pt]{geometry}
\usepackage[english]{babel}
\usepackage[cp1250]{inputenc}
\usepackage[T1]{fontenc}
\usepackage{booktabs}
\usepackage[title]{appendix}
\usepackage{url}
\usepackage[maxbibnames=9]{biblatex}
\usepackage[pdftex]{graphicx}
\usepackage{graphicx}
\usepackage[usenames,dvipsnames]{xcolor}
\usepackage{amsmath}
\usepackage{fixmath}
\usepackage{sidecap}
\usepackage{wrapfig}
\usepackage{subfig}
\usepackage{pdfpages}
\usepackage{hyperref}

\newtheorem{thm}{Theorem}

\newtheorem{rem}{Remark}

\title{Stability of fixed points in an approximate solution of the spring-mass running model}
\author{Zofia Wr\'oblewska\footnote{Corresponding Author: zofia.wroblewska@pwr.edu.pl}, Piotr Kowalczyk, {\L}ukasz P{\l}ociniczak\\\small{Faculty of Pure and Applied Mathematics},\\ \small{Wroclaw Universiry of Science and technology, Wroclaw, Poland}}

\begin{document}
\maketitle
\begin{abstract}	
We consider a classical spring-mass model of human running which is built upon an inverted elastic pendulum. Based on our previous results concerning asymptotic solutions for large spring constant (or small angle of attack), we construct analytical approximations of solutions in the considered model. The model itself consists of two sets of differential equations - one set describes the motion of the centre of mass of a runner in contact with the ground (support phase), and the second set describes the phase of no contact with the ground (flight phase). By appropriately concatenating asymptotic solutions for the two phases we are able to reduce the dynamics to a one-dimensional apex to apex return map. We find sufficient conditions for this map to have a unique stable fixed point. By numerical continuation of fixed points with respect to energy, we find a transcritical bifurcation in our model system. \\
		
		
\end{abstract}

\section{Introduction}
\label{indroduction}
\noindent
Running is one of the most fundamental means of bipedal animal locomotion. Superficially, it may seem mundane. However, it is a result of complex interactions between neural and muscular systems \cite{Dic10, Gor17}. Remarkable efficiency and stability of animal running over rough terrain attracts the attention of scientists throughout various disciplines spanning from biomechanics, medicine, applied mathematics, to robotics \cite{Bie18, Col05, Hol06}. Running can also be considered from several other points of view that are not of biomechanical nature. For instance, in competitive sports, one is usually interested in optimization problems based on traversing a given distance in the shortest time. A classical model of Keller \cite{Kel73}, which was based on early research of Hill \cite{Hil25}, provides an interesting account of this optimization. Some modern approaches to that subject can be found in \cite{Woo91, Pri93, Aft19}.  

To understand the very process of running, one usually constructs a model that can capture its relevant properties. This can be done on a spectrum of various complexity levels. Here, we are concerned with the, so-called, spring-mass running model for bouncing gaits that has been introduced by Blickhan \cite{blick} and thoroughly investigated by MacMahon and Cheng \cite{mcmahon}. This model can be considered conceptual in the sense that it is a simple realisation of the very motion of a hopping leg. It consists of a spring-loaded inverted pendulum (SLIP) that swings during the phase of contact with the ground (support phase) and is then  ejected upwards into the aerial phase (flight phase). The exceptional accuracy of that model is based on the validity of the reductionist approach. It stems from dimensional analysis of fundamental physical quantities describing the gait of various, not necessarily bipedal, animals. Some experimental studies concerning the SLIP model have been conducted for example in \cite{Far91, Far93, He91}. This model, due to its robustness, has also been extensively used in robotics to design and construct various legged robots (see the seminal paper by Raibert \cite{Rai86} and a modern review \cite{Agu16}). Furthermore, in the literature one can find a number of important generalizations of the SLIP model such as varying attack angle \cite{Gan18}, introducing control \cite{Sat04, Tak17, Sha16}, damping \cite{Sar10}, and multi-legged hoppers \cite{Gan14}. 

In spite of being a relatively simple model, SLIP does not posses a closed form analytical solution. Therefore, in order to study its properties one usually uses quantitative dynamical systems analysis or numerical methods. To gain more insight, an approximate analytical solutions can be obtained. This has been done in the work of Geyer, Seifarth, and Blickhan \cite{geyer} which is one of the basis of our subsequent results. Moreover, some approaches to analysis of the current model based on approximate solutions have been conducted in \cite{Sch00, Ghi05, Sar10} and in our earlier works \cite{plocin, okras}. There, we have conducted a rigorous Poicar\'e-Lindstedt asymptotic analysis based on the perturbation theory in the regime of large spring stiffness and small angle of attack. Moreover, we have also found some expansions in the limit of large horizontal velocity. A thorough review of mathematical studies of various locomotion models has been conducted in \cite{Hol06}, to which we refer the interested reader. 

In this paper, we obtain analytical solutions to the model by concatenating our approximate solutions in the support phase with the solutions of the flight phase, with the latter phase representing parabolic trajectories. Using these solutions we obtain a one-dimensional Poincar\'e map that takes an apex in the flight into the next apex. We then investigate the existence of fixed points of the Poincar\'e map so constructed and find sufficient conditions for the stability of fixed points. This shows that this model can describe a stable running gait. Our approach is similar to the one conducted in \cite{geyer} where the authors use a different approximation strategy stemming from the principles of mechanics, which, however, obscures the error one makes by necessarily approximating the solutions. In contrast to that, our approach is based on rigorous asymptotic analysis that clearly formulates the region of validity of the model \cite{plocin,okras}. Moreover, as our numerical studies indicate, the fixed point of the Poincar\'e map undergoes a transcritical bifurcation with increasing energy. A thorough investigation of this phenomenon, as well as the investigations of the existence of fixed points corresponding to asymmetric solutions, will be the subject of our future studies. 

The rest of the paper is outlined as follows. In Sec.~\ref{massmodel} we introduce the governing equations of the spring-mass model of running. Then, in Sec.~\ref{approx}, we introduce approximate solutions of our model equations obtained by means the Poincar\'e - Lindstedt method. These equations are then used in Sec.~\ref{stability} to obtain a one-dimensional map which is also analyzed in this section. Finally, Sec.~\ref{conclusion} concludes the paper.   
\section{Spring-mass running model}
\label{massmodel}

\begin{figure}
\center
\includegraphics[width=0.7\textwidth]{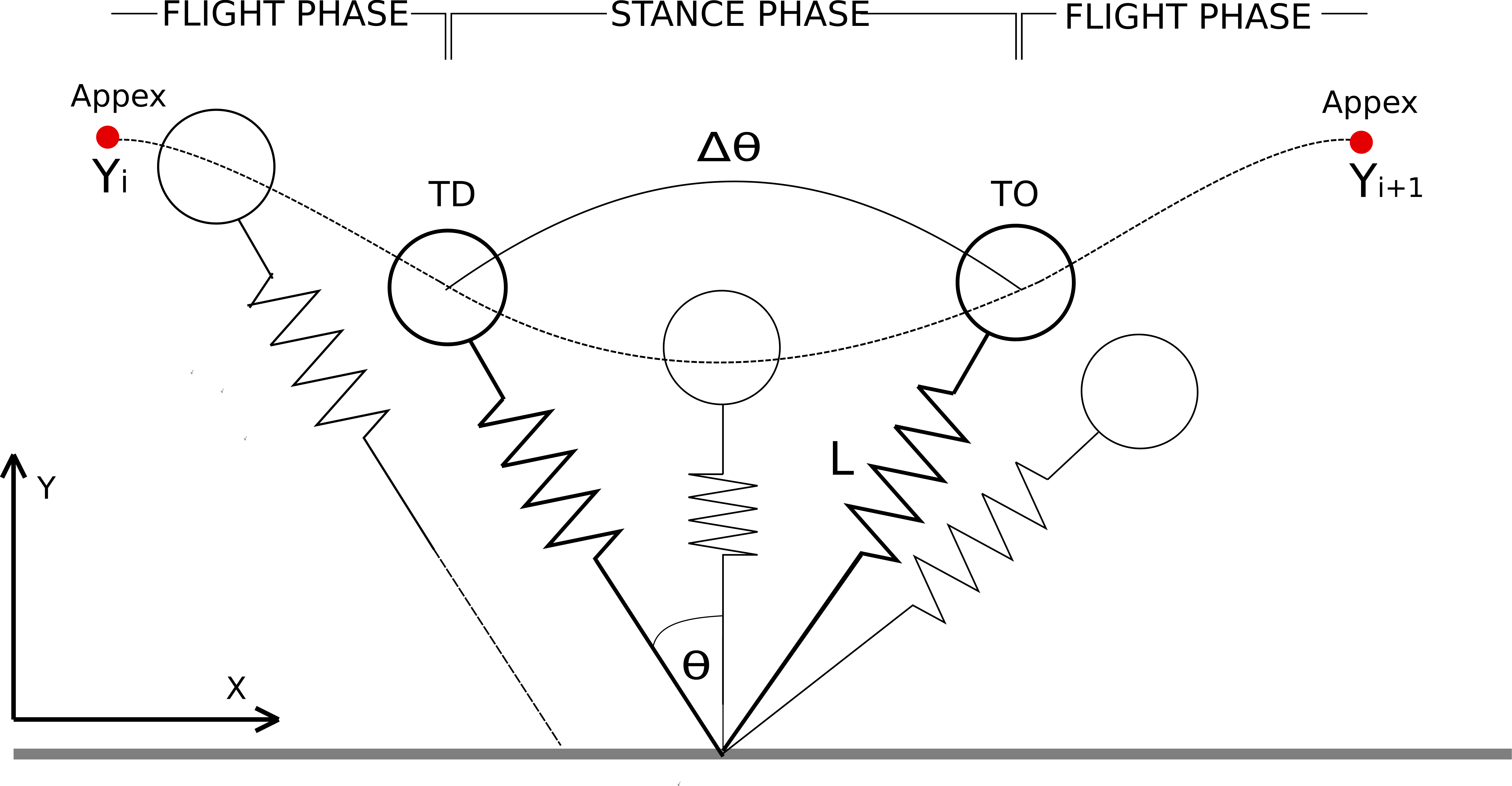}

\caption{\small Schematic of the spring-mass model for running. Dimensionless parameters: $(X,Y)$ -  Cartesian coordinates of the point mass, $(L,\theta)$ - polar coordinates, where $L=\sqrt{X^2+Y^2}$ is the radial while $\theta$ is the angular position of the point of mass, and additionally $\Delta\theta$ is angle swept during stance. Moreover, we assume that $\theta_{TD}=-\alpha$, where $\alpha \in (0,\pi/2)$.}
\label{figure1}
\end{figure}

\noindent
In this section we will derive  nonlinear equations, which describe running during the stance phase (see Fig. \ref{figure1} for notation). 
Using  nondimensional polar coordinates $(L,\theta)$, the Lagrange function of
the stance phase is given by 
\begin{equation}
	{\cal L}=\frac{1}{2}\left( (L')^2 + L^2(\theta')^2\right) - \frac{K}{2}\left( 1-L\right)^2 -L\cos\theta,
	\label{lagrange}
\end{equation}
where the  parameter $K$  denotes the nondimensional spring stiffness. So the two Euler-Lagrange equations are as follows
\begin{equation}
    \left\{
	\begin{array}{lcl}
	\dfrac{d}{dt}\left( \dfrac{\partial{\cal L}}{\partial L'}\right)=\dfrac{\partial{\cal L}}{\partial L} & \Longrightarrow & 
	L''=L (\theta')^2 + K (1-L) -\cos\theta,\vspace{6pt}\\
	\dfrac{d}{dt}\left( \dfrac{\partial{\cal L}}{\partial\theta'}\right)=\dfrac{\partial{\cal L}}{\partial \theta} & \Longrightarrow & 
	\left( L^2\theta'\right)'= P' = L\sin\theta,\\
	\end{array}
	\right.
	\label{twoEL}
\end{equation}
where  $P=L^2\theta'$  denotes the  nondimensional angular momentum.  

 Observe that the dimensionless form of the mechanical energy during the contact phase denoted by $E_s$  (see (\ref{lagrange})) is given by
\begin{equation}
	E_s=\frac{1}{2}\left[ (L')^2 + L^2(\theta')^2 + K(1-L)^2 \right ] + L\cos\theta.
	\label{energyd}
\end{equation}
On the other hand, the total energy $E_s$ given by  equation (\ref{energyd}), for small angles $\theta$, takes an alternative form
\begin{equation}
	\widetilde{E}_s=\frac{1}{2}\left[ (L')^2 + L^2(\theta')^2 + K(1-L)^2 \right ] + L.
	\label{energyalter}
\end{equation}
Consequently, the Lagrangian (\ref{lagrange}) is also modified. We should note that for small angles approximation, change of the angular momentum of the point mass we consider is zero, that is 
\begin{equation}
\frac{d}{dt} P = \frac{d}{dt}\left( L^2\theta'\right)= 0.
	\label{momentum}
\end{equation}
From equations (\ref{energyalter}) and (\ref{momentum}), we now obtain a set of conditions that relate the system state at take-off to the state  at touch-down. During take-off and touch-down, denoted by {\it TO} and {\it TD}, respectively, we have the radial symmetry with  rest length $L=1$ at each phase transition
\begin{equation}
	L_{TO}=L_{TD}=1,
	\label{stab1}
\end{equation}
and so from  (\ref{energyalter}) we have then
\begin{equation}
	2\left(\widetilde{E}_s-1\right)= \left(L'\right)^2_{TO/TD} + \left(\theta'\right)^2_{TO/TD}.
	\label{energyalter2}
\end{equation}
Furthermore, after integration  (\ref{momentum}), the following equation  $\; \theta'_{TO/TD}=cL^{-2}_{TO/TD}=c$ holds, where $c$ is some real constant. In our case: forward locomotion the constant $c$ is positive.
Thus 
\begin{equation}
    \theta'_{TO}= \theta'_{TD}>0,
	\label{stab2}
\end{equation}
and $\;2\left(\widetilde{E}_s-1\right) - c^2= \left(L'\right)^2_{TO/TD},\;$ and hence $|L'_{TO}|=|L'_{TD}|$. We can now write
\begin{equation}
	L'_{TO}=-L'_{TD}>0.
	\label{stab3}
\end{equation}
Note that in the full system, in addition to (\ref{stab1}), (\ref{stab2}) and (\ref{stab3}), we require (see also Fig. \ref{figure1})
\begin{equation}
	\theta_{TO}=  \Delta \theta + \theta_{TD}, \quad (\theta_{TD}<0, \; \theta_{TO}>0),
	\label{stab4}
\end{equation}
where only the angle $\Delta \theta$ swept in stance cannot
simply be expressed by the state at the time of touch-down. Moreover, it should be observed that for small  angles $\theta$ the conservation of energy and angular momentum do not enforce symmetry in the stance phase i.e. through the arbitrary value of $\theta$ parameter, asymmetric solutions are also allowed. 

In the next section we will derive the formula for $\Delta \theta$.

\section{Constructing approximate solutions}
\label{approx}
\noindent
 An asymptotic analysis of the main equations (\ref{twoEL}) with the use of the Poincar\'e - Lindstedt series was carried out in \cite{plocin}.  The solution of the problem is based on
the perturbative expansion related to the significant spring stiffness ($K \rightarrow \infty$). To simplify matters we set $\epsilon=1/\sqrt{K}$ ($\epsilon \rightarrow 0^{+}$) and introduce \begin{equation}
  \begin{array}{c}
 \widetilde\omega(\epsilon)=1-\frac{1}{2}\epsilon^2\theta_d^2.   
  \end{array}
\end{equation}
In what follows we will operate on the order of $O(\epsilon^2)$ as $\epsilon\rightarrow 0^+$ truncating higher order terms. All equations should then be understood as approximations in this asymptotic sense. For clarity we will retain using the equality sign instead of more rigorous approximate equality after performing the truncation. 

\begin{itemize}
\item Radial motion during stance $L(t)$ (see (23) in \cite{plocin}):\\
The $L(t)$ approximation is as follows
 \begin{equation}
	L(t)= 1-\epsilon L_d\sin{\left(\epsilon^{-1} \widetilde\omega(\epsilon)t\right)} + \epsilon^2\left(\theta_d^2-\cos{\alpha}\right)\left(1-\cos{\left(\epsilon^{-1} \widetilde\omega(\epsilon)t\right)}\right),
	\label{radial55}
\end{equation}
where the initial conditions  have the form for $\alpha>0$
\begin{equation}
  \begin{array}{ll}
	\theta(0)=\theta_{TD}=-\alpha, & L(0)=L_{TD}=1, \\
	 \theta'(0)=\theta'_{TD}=\theta_d, & L'(0)=L'_{TD}=-L_d,
	 \label{initialc}
  \end{array}
 \end{equation} 
with
	\begin{equation}
	\label{ThetaLdUV}
	\theta_d = X' \cos\alpha- Y' \sin\alpha, \quad L_d = X' \sin\alpha + Y' \cos\alpha.
	\end{equation}
$X'$ and $Y'$ are the horizontal and vertical  velocities  made dimensionless by the factor $(\textrm{acceleration of gravity} \times  \textrm{length})^{1/2}$, which are called Froude numbers in fluid mechanics.  In addition, it is worth noting that from (\ref{energyd}) 
\begin{equation}
	 L_d=|L'_{TD}|=\sqrt{2E_s-\theta_d^2-2\cos{\alpha}},
	\label{ab1}
\end{equation}
where $\theta_d$ refers to the angular momentum $P$ at touch-down ($TD$).

The general solution (\ref{radial55}) of the radial motion $L(t)$ during stance describes a sinusoidal oscillation around $L=1+A(\epsilon)$ with amplitude $B(\epsilon)$ and frequency $\epsilon^{-1} \widetilde\omega(\epsilon)$, where
\begin{equation}
 A(\epsilon)=\epsilon^2\left(\theta_d^2-\cos{\alpha}\right) \quad \mbox{and} \quad B(\epsilon)=\epsilon\sqrt{\epsilon^{2}(\theta_d^2-\cos{\alpha})^2 + L_d^2}.
	\label{abformula3}
 \end{equation}
But, as  shown in Fig. \ref{Fig2} on the left, the solution only holds for $L(t)\leq 1$.

In formula  (\ref{radial55}), the parameter $t$ ranges from 0 to the contact time $t_C$ given by
\begin{equation}
	t_{C}=\frac{2\epsilon}{\widetilde\omega(\epsilon)}\arccos{\frac{A(\epsilon)}{B(\epsilon)}}=\frac{4\epsilon}{2-\epsilon^2\theta_{d}^2}\arccos{C(\epsilon)},
	\label{Tc22}
\end{equation}
where
\begin{equation}
 C(\epsilon)=\frac{A(\epsilon)}{B(\epsilon)}=\frac{\epsilon\left(\theta_d^2-\cos{\alpha}\right)}{\sqrt{\epsilon^{2}(\theta_d^2-\cos{\alpha})^2 + L_d^2}},
	\label{abformula2}
 \end{equation}
 and $\epsilon \rightarrow 0^{+}$.
Observe that $L(0)=L(t_C)=1$ and $L(t_C/2)=1-\Delta L_{MAX}$ (see Fig. \ref{Fig2}). The maximum leg deflection during the stance phase, denoted by  $\Delta L_{MAX}$, is given by the difference of the amplitude $B(\epsilon)$ of  the radial motion $L(t)$ and the shift $A(\epsilon)$ of the touch-down position, i.e.
\begin{equation}
\Delta L_{MAX}= B(\epsilon) - A(\epsilon).
	\label{maksspring}
 \end{equation}
So restriction  to small values
of $L-1$ is adequately formulated by $B(\epsilon)-A(\epsilon) = O(\epsilon) \ll 1$ as $\epsilon \rightarrow 0^{+}$. For illustrations, we refer to Fig. \ref{Fig2}. In the case of the data from  Fig. \ref{Fig2}, we get $A(\epsilon) = -0.0016$ and $B(\epsilon) = 0.0515$, where $\epsilon=1/\sqrt{15}$. Hence $\Delta L_{MAX}=0.0531$, and we can check, that $\Delta L_{MAX}=1-L(t_C/2)$, where $t_C=0.8558$.  In addition, note that the shift $A(\epsilon)$ may have
 also a positive result, when $\theta_d>\sqrt{\cos{\alpha}}$, then it will move $1$ below  $1+A(\epsilon)$ and reduce the maximum spring compression $\Delta L_{MAX}$.

\begin{figure}
		\centering
		\includegraphics[width=0.49\textwidth]{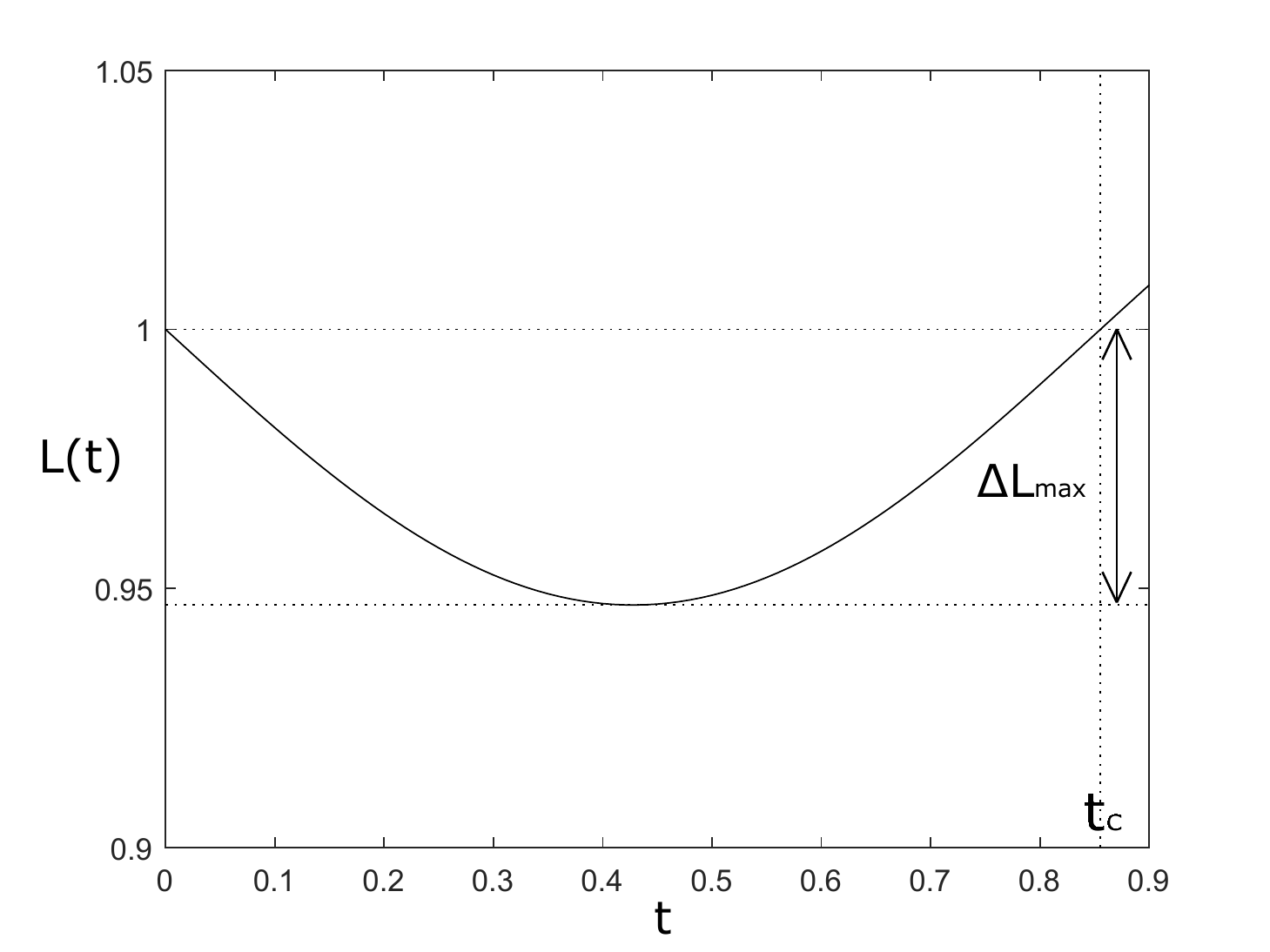}
		\includegraphics[width=0.49\textwidth]{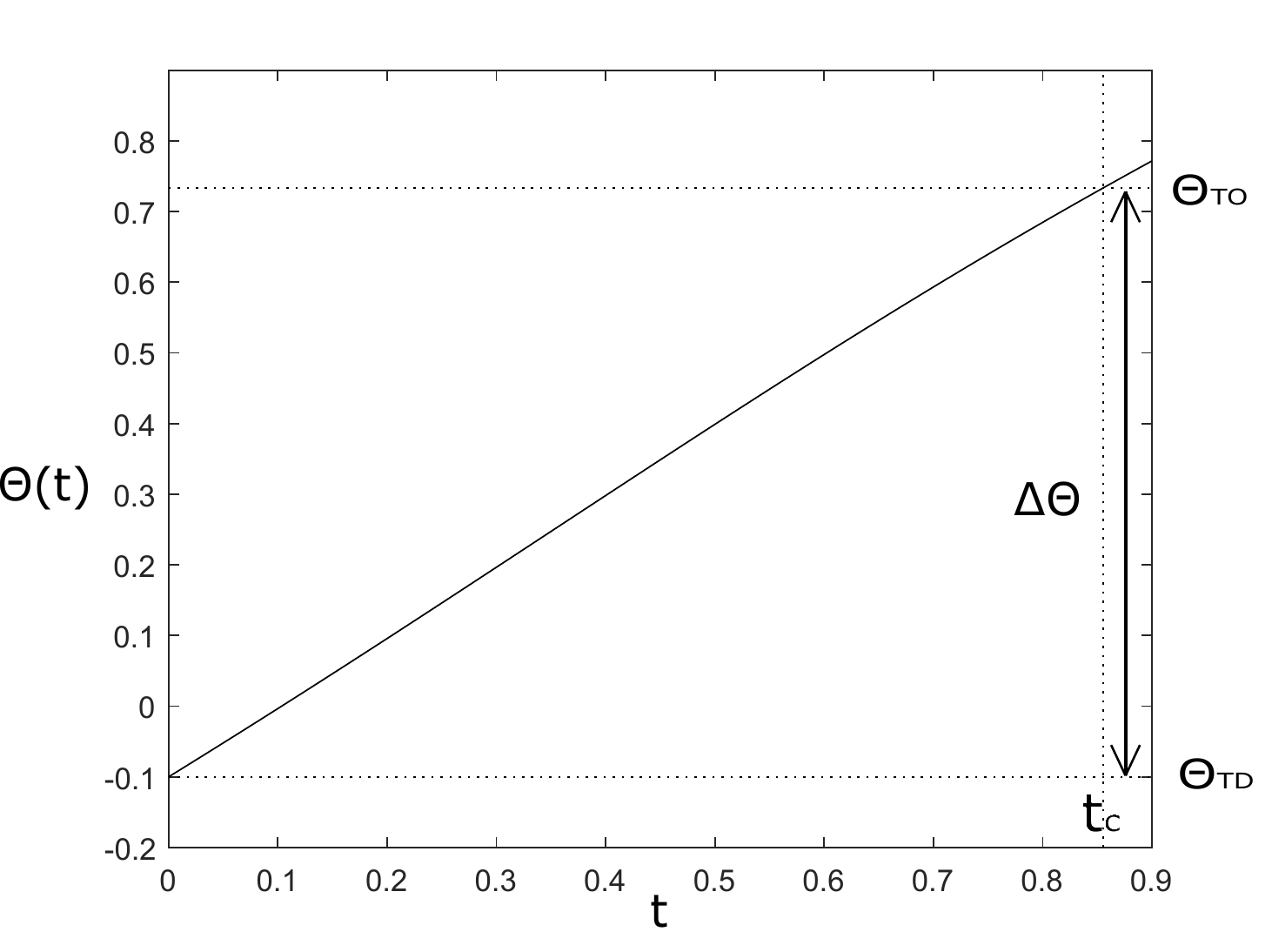}
		\caption{\small   $L$ and $\theta$ approximations  during the stance phase. Here $\alpha=0.1$, $X'=1$, $Y'=0.1$ and $K=15$.}
		\label{Fig2}
\end{figure}

\item Angular motion during stance $\theta(t)$ (see (24) in \cite{plocin}):\\
The $\theta(t)$ approximation with  the initial conditions (\ref{initialc}) is as follows
 \begin{equation}
	\theta(t)= - \alpha +\theta_d\widetilde\omega(\epsilon)t-\frac{1}{2}\left(\widetilde\omega(\epsilon)t\right)^2\sin{\alpha} +2\epsilon^2L_d\theta_d\left(1-\cos{\left(\epsilon^{-1} \widetilde\omega(\epsilon) t\right)}\right).
	\label{angular22}
\end{equation} 
Observe that $\theta(0)=\theta_{TD}=- \alpha$ and $\theta(t_C)=\theta_{TO}=\Delta \theta +\theta_{TD} = \Delta \theta - \alpha$  (see Fig.  \ref{Fig2} on the right).

\item Angle swept during stance:\\
From formulas (\ref{angular22}) and (\ref{Tc22}) we obtain the angle $\Delta \theta$ swept during stance
\begin{equation}
		\Delta \theta= 	\theta(t_C) +\alpha =2\epsilon\theta_d \times\left[\arccos{C(\epsilon)} + 2\epsilon L_d\left(1-C^2(\epsilon)\right)\right]
		- 2\epsilon^2\sin{\alpha}\left(\arccos{C(\epsilon)}\right)^2,
	\label{angular222}
\end{equation} 
where $C(\epsilon)$ is already defined in (\ref{abformula2}).
Therefore, the angle swept $\Delta \theta$ is
uniquely defined by the system state at touch-down $(L'_{TD}, \theta_{TD}, \theta'_{TD})$ and the spring stiffness parameter  $K$. Moreover, the landing angle
$\theta_{TD}=- \alpha$ affects $\Delta \theta$ also  by determining  the radial and angular landing velocity.   

For the data in Fig. \ref{Fig2} we get $\theta(t_C)= \theta_{TO} = 0.7334$. Hence $\Delta \theta=\theta_{TO} - \theta_{TD}=0.8334.$ Which agrees with the value calculated in formula (\ref{angular222}).

Finally, let's observe, that 
the Taylor expansion of (\ref{angular222}) to the second order, about $\epsilon=0$, is given by 
\begin{equation}
		\Delta \theta= \pi\theta_d \epsilon + \left[4\theta_d L_d - 2\frac{\theta_d}{L_d}\left(\theta_d^2-\cos{\alpha}\right) - \frac{1}{2}\pi^2\sin{\alpha}\right]\epsilon^2,
	\label{angular333}
\end{equation} 
where $\epsilon \rightarrow 0$. Since $\Delta \theta$ approximation  given by formula (\ref{angular222})  is of the order $\epsilon^2$, all terms of the order $O(\epsilon^3)$ in formula (\ref{angular333}) are omitted. 
\end{itemize}

An important case of stability analysis, considered  also by Geyer, Seyfarth and Blickhan (cf.  \cite{geyer}), is 
\begin{equation}
C(\epsilon) = 0 \quad \Longleftrightarrow \quad \theta_d=\sqrt{\cos{\alpha}} \quad \mbox{and} \quad L_d=\sqrt{2E_s-3\cos{\alpha}}.		
	\label{special}
\end{equation} 
In  the special case we get the shift $A(\epsilon)=0$ and  the leg deflection parabola is equal to the amplitude i.e. $\Delta L_{MAX} =B(\epsilon)$. Moreover, the angular velocity at touch-down is 
$\theta_d=\sqrt{\cos{\alpha}}$ and the contact time in formula (\ref{Tc22}) is equal to $2\pi\epsilon/(2-\epsilon^2\theta_{d}^2)$.

Putting $C(\epsilon) = 0$ into the equation (\ref{angular222}) or equivalently $\theta_d=\sqrt{\cos{\alpha}}$ into the equation (\ref{angular333}), we get (see also (\ref{ab1}))
\begin{equation}
		\Delta \theta= \pi\theta_d \epsilon + \left[4\theta_d L_d  - \frac{1}{2}\pi^2\sin{\alpha}\right]\epsilon^2 =\pi\sqrt{\frac{\cos{\alpha}}{K}} + \left[4\sqrt{\cos{\alpha}(2E_s-3\cos{\alpha})}  - \frac{1}{2}\pi^2\sin{\alpha}\right]\frac{1}{K}. 
	\label{angular666}
\end{equation} 
Parameter combinations $\{\alpha, E_s(\alpha), K(\alpha,E_s)\}$ leading to  a symmetrical stance phase, i.e. $\Delta \theta=2\pi \;$, will be given, for the case of $C(\epsilon) = 0$ in Sec.~\ref{a0}.

By setting the stance phase limits to
symmetry,  we derive in the next section approximations of the dimensionless leg stiffness $K$.

\subsection{Symmetric solutions}
\label{symmetric}
\noindent
 By limiting ourselves to the second order expansion (\ref{angular333}), the solution $\Delta \theta = 2\alpha$ comes down to the quadratic equation for $\sqrt{K}$. 
A physically reasonable solution gives the following expression to approximate the spring stiffness
\begin{equation}
\label{approx444}
K =\frac{\left[\pi\theta_d + \sqrt{\Delta}\right]^2}{16\alpha^2},
\end{equation}
where
\begin{equation}
\label{delta3}
\Delta= \pi^2\theta_d^2 + 8\alpha\left[4\theta_d L_d - 2\frac{\theta_d}{L_d}\left(\theta_d^2-\cos{\alpha}\right) - \frac{1}{2}\pi^2\sin{\alpha}\right] .
 \end{equation}
 Using  $\alpha  \rightarrow 0^+$ in  (\ref{delta3}), the second component tends to zero and we get \footnote{Since the angular velocity  at touch-down $\theta_d(\alpha)$ may depend on the angle $\alpha$, to be asymptotically consistent we should write
$$
\widetilde{K} = \left(\frac{\pi \theta_d(0)}{2\alpha}\right)^2 \quad \text{as} \quad \alpha\rightarrow 0^+,$$ where $\theta_d(0)=X'$. For simplicity we use the formula (\ref{napprox555}).}

\begin{equation}
\label{napprox555}
\widetilde{K} = \left(\frac{\pi \theta_d}{2\alpha}\right)^2. 
\end{equation} 
The same approximation of $K$  as (\ref{napprox555})  was obtained in \cite{plocin} and \cite{wrobel}. Moreover, in \cite{okras} it was shown that the leading order of the expansion of $K$ as $\alpha  \rightarrow 0^+$ is indeed $(\pi X')^2/(2\alpha)^2$.

In the symmetric case we have a particular boundary value problem to solve. Let $(\theta(t,K), L(t,K))$ be the solution of the system (\ref{twoEL}) with (\ref{initialc}). Find $K^*$ and the smallest time $t^*>0$ satisfying 
\begin{equation}
\label{brzeg}
\theta(t^*,K^*)=\alpha, \quad\quad L(t^*, K^*) = 1.
\end{equation}
The above problem can be easily solved numerically using the shooting method, as was done for example in \cite{mcmahon}. Graphical illustrations confirming the validity (\ref{napprox555}) as an approximation of the solution to the boundary problem (\ref{brzeg}) are presented in \cite{plocin} and \cite{wrobel}.

\begin{figure}
		\centering 
		\includegraphics[width=0.49\textwidth]{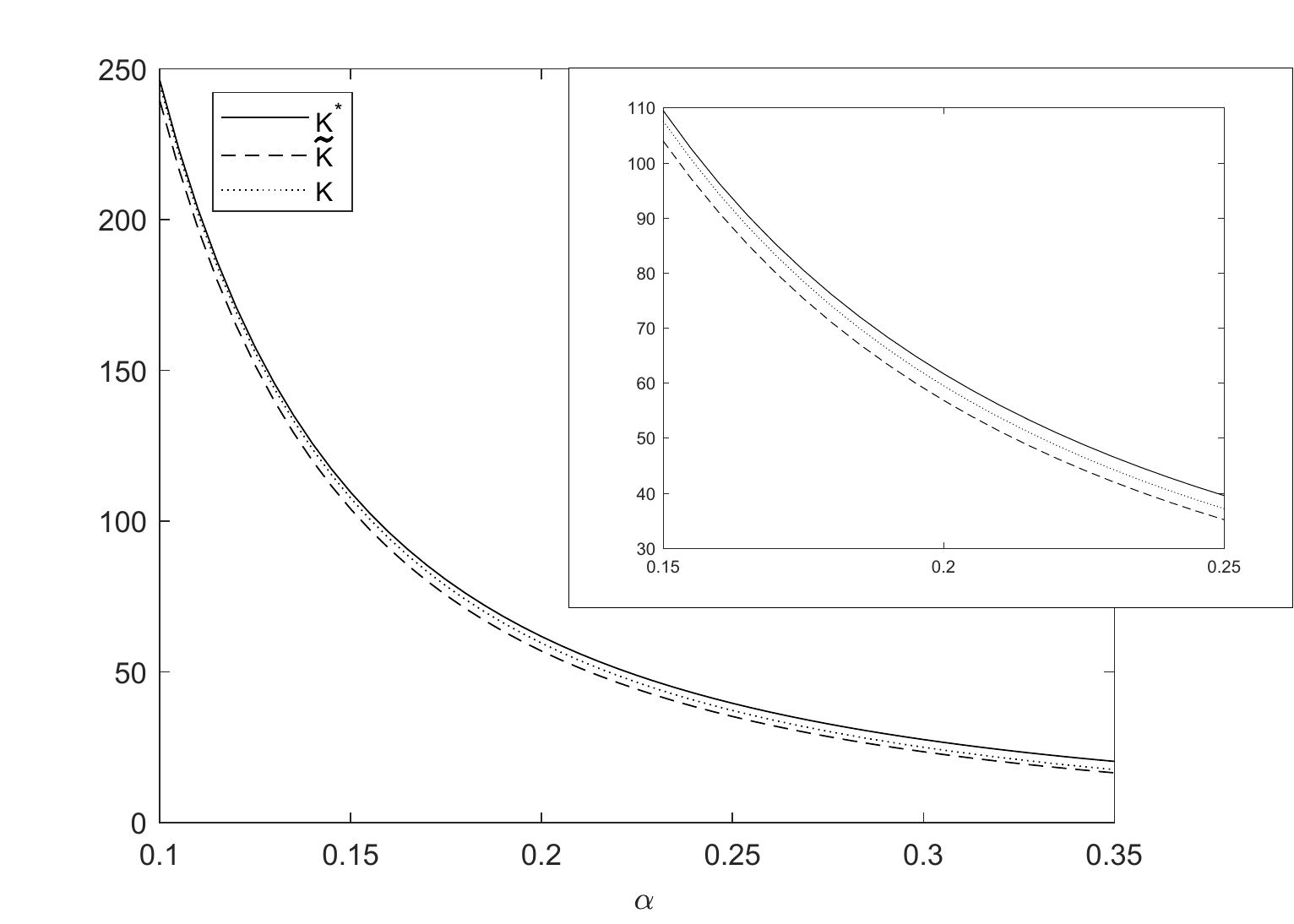}
		\includegraphics[width=0.49\textwidth]{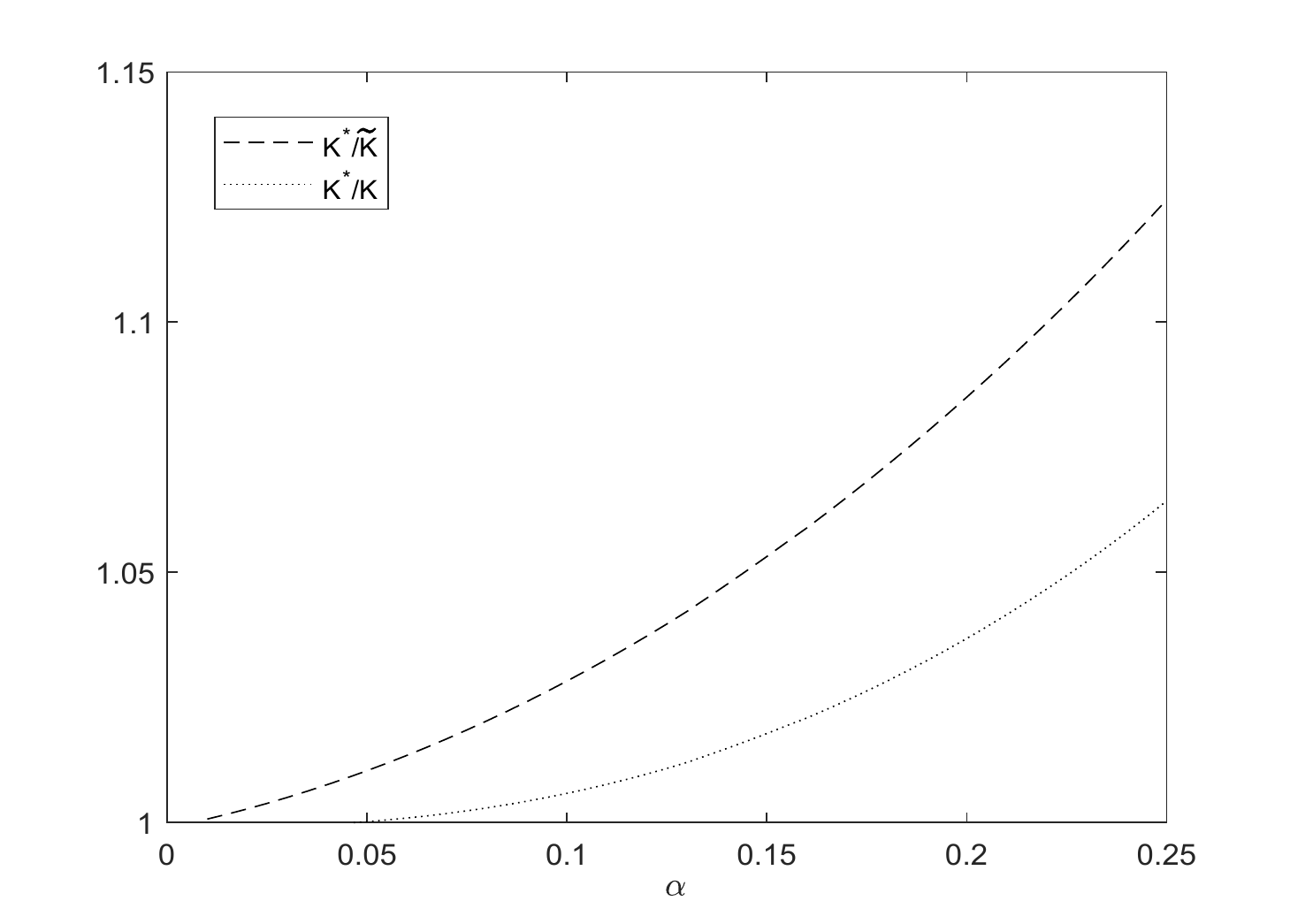}
		\caption{\small On the left: comparison  between the numerically calculated value of $K^*$ and its $K$ approximation, calculated from (\ref{approx444}), and the second $\widetilde{K}$ approximation, calculated from (\ref{napprox555}), for varying $\alpha$ from 0.1 to 0.35. On the right: ratio of $K^*$ 
to $K$ and $\widetilde{K}$ for varying $\alpha$ from 0.01 to 0.25. Here $X'=1$ and $Y'=0.1$.}
		\label{Fig3}
\end{figure}

In Fig. \ref{Fig3}, on the left, we can see what the approximations $K$ given by (\ref{approx444}) and $\widetilde{K}$ given by (\ref{napprox555}) look like. Obviously, if $\alpha\rightarrow 0^+$  then $K^*\rightarrow \infty$ and both approximations of $K^*$ are very good. However, for moderate values of parameter $\alpha$, $K$ is slightly better than $\widetilde{K}$.
Moreover the right side of Fig. \ref{Fig3} illustrates the ratio of $K^*$  to both approximations $K$ and $\widetilde{K}$  for varying $\alpha$. Indeed, the approximation $K$ behaves better than $\widetilde{K}$ in the presented range of the parameter $\alpha$, i.e. $0.01 \leq \alpha \leq 0.25$. It turns out that the approximation $K$ works very good for $\alpha<0.05$ while $\widetilde{K}$ is accurate only for $\alpha<0.01$.

\section{Stability of spring-mass running}
\label{stability}

\subsection{Analytical apex return map}
\label{apex}
\noindent
At the beginning of this section we will deduce  the analytical
solution to calculate dependencies between
two  subsequent apex heights.
Identification of the system during the phase transitions: from flight to stance and from stance to flight will allow the derivation of the  apex return map $Y_{i+1}(Y_i)$.
In Sec.~\ref{massmodel} we showed that assuming small angles, the system state at take-off ($TO$) is related to the state
at touch-down (TD) through a set of the equations  (see (\ref{stab1}), (\ref{stab2}), (\ref{stab3}) and (\ref{stab4}))
\begin{equation}
\begin{array}{l}
	L_{TO}=L_{TD}=1,\\
	L'_{TO}=-L'_{TD}=L_d,\\
	\theta_{TO}= \Delta \theta +\theta_{TD}= \Delta \theta - \alpha,\\
	\theta'_{TO}= \theta'_{TD}=\theta_d.
\end{array}	
	\label{stabw2}
\end{equation}
Starting with the correct initial values, the mapping between the height $Y_i$ of the apex $i$ and the touch-down state in dimensionless polar coordinates is expressed as
\begin{equation}
\begin{array}{lclcl}
  Y_i & \longrightarrow & 
 \left [
 \begin{array}{l}
	X'=\sqrt{2(E_s-Y_i)}\\
	Y=\cos{\alpha} \\
	Y'=-\sqrt{2(Y_i-Y)}
\end{array}	
 \right]_{TD}
 & \longrightarrow & 
 \left [
\begin{array}{l} 
	L=1\\
	L'=X'\sin{\theta} + Y'\cos{\theta}\\
	\theta =-\alpha\\
	\theta'=X'\cos{\theta} - Y'\sin{\theta}
\end{array}	
 \right]_{TD},
\end{array}	
	\label{map}
\end{equation}
where $E_s$ is the corresponding dimensionless energy of the system prior to touch-down, given from (\ref{energyd}) by $E_s = (\theta^2_{d} + L^2_{d})/2 +\cos{\alpha}$. Since the Froude numbers $X'$ and $Y'$ depend on $Y_i$, then the radial velocity $L_d$ and the angular velocity $\theta_d$ also depend on $Y_i$ (see (\ref{ThetaLdUV}) or (\ref{map})). 
Then, an additional mapping is required between the take-off state and the height $Y_{i + 1}$ of the apex $i + 1$, which is
\begin{equation}
\begin{array}{lcl}
 \left [
 \begin{array}{l}
	X'=L'\sin{\theta} + \theta'\cos{\theta}\\
	Y=\cos{\theta} \\
	Y'=L'\cos{\theta} - \theta'\sin{\theta}
\end{array}	
 \right]_{TO}
 & \longrightarrow & 
 \left [
\begin{array}{l}
	X'_{i+1}=X'_{TO}\\
    Y_{i+1}=Y_{TO}+\frac{1}{2}\left ( Y'_{TO}\right )^2
\end{array}	
 \right].
\end{array}	
	\label{map2}
\end{equation}
By using both mappings (\ref{map}) and (\ref{map2}), the apex return map function $Y_{i+1}(Y_i)$
can be constructed
and takes the following form after some simplifications
\begin{equation}
\begin{array}{l}
	Y_{i+1}(Y_i) =\cos{(\alpha - \Delta\theta_i)} + \left [ \sin{(2\alpha - \Delta\theta_i)}\sqrt{E_s-Y_i} +
	\cos{(2\alpha - \Delta\theta_i)}\sqrt{Y_i- \cos{\alpha}}  \right ]^2, 
\end{array}		
	\label{map3}
\end{equation}
where $\Delta\theta_i = \Delta\theta(\theta_d,L_d)$ is  given by (\ref{angular222}) or (\ref{angular333}). To mark the periodicity for this running model, we introduced the index $i$ to the notation $\Delta\theta_i$.
So $\Delta\theta_i$
is the angle swept during stance between two subsequent
flight phases, when  the apex height $Y_i$ transforms into $Y_{i + 1}$. 
Now you can see how $\Delta\theta_i= \Delta\theta(\theta_d,L_d)$  depends on $Y_i$.  From formulas (\ref{ThetaLdUV}) and (\ref{map}) we have 
\begin{equation}
	\label{ThetaLdUV2}
	\begin{array}{r}
	\theta_{d} =\sqrt{2}\left [ \cos\alpha\sqrt{E_s-Y_i}\right.\\ \left.-\sin\alpha\sqrt{Y_i-\cos{\alpha}} \right ]
\end{array}, \quad 
\begin{array}{r}
	L_{d} = \sqrt{2}\left [ \sin\alpha\sqrt{E_s-Y_i}\right.\\ \left. + \cos\alpha\sqrt{Y_i-\cos{\alpha}} \right ]
\end{array}.
\end{equation}
Moreover, observe that (see (\ref{ab1}))
\begin{equation}
\theta^2_{d} + L^2_{d}=2(E_s-\cos{\alpha}).
	\label{map3333}
\end{equation}
If the apex height $Y_{i+1}$ is less then landing height ($Y_{i+1} < \cos{\alpha}$), the leg will get into the ground, and the stumble will occur. So the apex return map only exists  otherwise. 

Regardless of the angle swept during stance, from equations (\ref{stabw2}) the kinetic energy $\frac{1}{2}\left( (L')^2 + L^2(\theta')^2\right)$ is equal at touch-down and take-off.
For asymmetric contact phases, since $Y_{TD} \neq Y_{TO}$, there is a net change in system energy $\Delta E=Y_{TO}-Y_{TD}$. On the other hand,  during symmetric stance phases ($Y_{TD}=Y_{TO}$) the corresponding shifts in energy at touch-down and take-off compensate each other, so  the system energy $E_s$ is constant. To restore the conservative nature of the model, change resulting from asymmetry, needs to be corrected. Thus, the appropriate horizontal velocity in (\ref{map2})
\begin{equation}
	X'_{i+1}=\sqrt{2(E_s-Y_{i+1})}. 
	\label{velocity}
\end{equation}
must be given during take-off so that the kinetic energy ($E_s-Y_{i+1}$) reduces the lack of potential energy at the end of stance phase. For the new apex hight $Y_{i+1}$ from (\ref{map3}) the adjusted $E_s$ is used.

\subsection{Existence of fixed points}
\label{fixedpoints1}
In the following section we need to check whether apex states $Y_{i+1}$ restricted by symmetric contracts can be found.
Solution (\ref{map2}) with condition (\ref{stabw2}) for  $\theta_ {TO} = \alpha$ leads to the corresponding apex height of  fixed points
\begin{equation}
	Y^*  = \cos{\alpha}+ \frac{1}{2}\left(\theta^*_{d}\sin{\alpha} - L^*_{d}\cos{\alpha}\right)^2,
\label{apexheight3}
 \end{equation}
where $\theta^*_d$ and $L^*_d$ represent the values of $\theta_d$ and $L_d$ given by (\ref{ThetaLdUV2}) at the fixed point $Y^*$.
Since   $\theta^*_d$ and $L^*_d$ are linked by the relationship (\ref{map3333}), the apex height of  fixed points can be expressed in the following equivalent form
\begin{equation}
	Y^*  = E_s - \frac{1}{2}\left(\theta^*_{d}\cos{\alpha} + L^*_{d}\sin{\alpha}\right)^2.
\label{apexheight3nowy}
 \end{equation}
 Substituting (\ref{apexheight3}) and (\ref{apexheight3nowy}) into equations (\ref{ThetaLdUV2}) yields
 \begin{equation}
\theta^*_{d} = \cos{\alpha}\left|\theta^*_{d} \cos{\alpha}+L^*_{d}\sin{\alpha}\right|\\
-\sin{\alpha}\left|\theta^*_{d}\sin{\alpha}-L^*_{d}\cos{\alpha}\right|,
\label{theta22}
\end{equation}
and also
\begin{equation}
L^*_{d} = \sin{\alpha}\left|\theta^*_{d} \cos{\alpha}+L^*_{d}\sin{\alpha}\right|\\
+\cos{\alpha}\left|\theta^*_{d}\sin{\alpha}-L^*_{d}\cos{\alpha}\right|.
\label{L22}
\end{equation}
It is easy to see that 
$$\theta^*_{d} \cos{\alpha}+L^*_{d}\sin{\alpha}>0.$$
So, it follows that equations (\ref{theta22}) and (\ref{L22}) hold
if
 \begin{equation}
\theta^*_{d}
\sin{\alpha} -  L^*_{d}\cos{\alpha} \le 0.
\label{22}
\end{equation}
From (\ref{ab1}) and (\ref{22}) we have
\begin{equation}
\theta^*_d \tan{\alpha}\le \sqrt{2E_s-(\theta^*_d)^2-2\cos{\alpha}}.
\label{Emin}
\end{equation}
Hence the system energy $E_s$ must satisfy the inequality
  $E_s \ge E_s^{min}$  with
\begin{equation}
E_s^{min}=\frac{(\theta^*_d)^2}{2\cos^2{\alpha}}+\cos{\alpha}.
\label{energymin}
 \end{equation}
 For the system energy $E_s = E_s^{min}$, the apex height of  fixed point $Y^*$ is $\cos{\alpha}$, which is equal to the landing height. 
 As the energy increases ($E_s > E_s^{min}$), the apex height increases, but it never exceeds the upper  boundary, i.e. $E_s$.

If the attack angle $\alpha$ is fixed, then the stiffness approximation given by  (\ref{approx444}) and (\ref{delta3}) is the lowest for the minimum energy system given by  (\ref{energymin}), i.e. 
\begin{equation}
K^{min} = \frac{\left[\pi\theta^*_d + \sqrt{\Delta_{min}
}\right]^2}{16\alpha^2},
\label{stiffnessaprox}
 \end{equation}
 where
 \begin{equation}
\Delta_{min}= \pi^2\left(\left(\theta^*_d\right)^2-4\alpha\sin{\alpha}\right) + 32\alpha\left(\theta^*_d\right)^2\tan{\alpha}- 16\alpha\cot{\alpha}\left(\left(\theta^*_d\right)^2-\cos{\alpha}\right).
 \label{delta5}
 \end{equation}
 
Whence, we proved that there is a fixed point under a certain assumption:
 \begin{thm} \textnormal{(On the existence of symmetric solutions)}. 
 Let the parameters $\alpha$, $E_s$  and $\theta^*_d$ satisfy the inequality
 \begin{equation}
E_s>\frac{(\theta^*_d)^2}{2\cos^2{\alpha}}+\cos{\alpha}, 
\label{energyminnowy}
 \end{equation}
 where $\alpha \in (0,\pi/2)$. Then there exist initial conditions for the stance phase 
 $$
 \left(\theta, \theta', L, L'\right) = \left(-\alpha, \theta_d^*, 1, \sqrt{2E_s-(\theta^*_d)^2-2\cos{\alpha}}\right)
 $$
 such that the fixed point equation (\ref{apexheight3}) and (\ref{apexheight3nowy}) are satisfied.
 \label{th1}
\end{thm}
 \begin{rem}
When the angular velocity at touch-down for the fixed point $Y^*$ equals $\theta_d^*=\sqrt{\cos\alpha}$, then the condition (\ref{energyminnowy}) from Theorem \ref{th1} takes the form 
 \begin{equation}
E_s > E_s^{min}= \frac{1}{2\cos{\alpha}}+\cos{\alpha}.
\label{energy222}
 \end{equation}
 Additionally we get
 \begin{equation}
K >  K^{min}= \frac{\left[\pi\sqrt{\cos{\alpha}} + \sqrt{\Delta_{min}
}\right]^2}{16\alpha^2},
\label{stiffness222}
 \end{equation}
 where
 \begin{equation}
\Delta_{min}= \pi^2\left(\cos{\alpha}-4\alpha\sin{\alpha}\right) + 32\alpha\sin{\alpha}.
 \label{deltamin222}
 \end{equation}
 \end{rem} 
\noindent The fixed points can only exist if a minimum energy $E_s^{min}$ and a minimum stiffness $K^{min}$ are
exceeded.

Considering $\theta_d^*=\sqrt{\cos\alpha}$ case, if $\alpha  \rightarrow 0$, the minimal system energy $E_s^{min}$ is approaching the value of $E_s=1.5$. 
To get a dimensional overview, the system energy then has the value: $mgl_0 E_s \approx 1125\: J$, where body mass $m = 75\: kg$, gravitational acceleration $g = 9.81\:  m/s^2$, and leg length $l_0= 1\:  m$. 
If additionally the initial  apex height  is set to $Y_0 = 1$, then in this case it corresponds to the initial velocity $x'=X'\sqrt{gl_0}=\sqrt{2(E_s-Y_0)gl_0}$, which is $3.13 \:  m/s$.

\subsection{Stability of fixed points}
\label{fixedpoints22}
Let us denote mapping (\ref{map3}) as $Y_{i+1}=f(Y_i)$.
Stable solutions of fixed point $Y^*$ fulfill the following condition 
\begin{equation}
|f'(Y^*)| < 1.
\label{stability2n}
 \end{equation}
 So to prove stability, at least one
 the parameter set $(\alpha, E_s, \theta^*_d)$  leading to solutions of $Y^*$ satisfying (\ref{stability2n}) should be identified.
 Starting with  (\ref{map3}), we get
\begin{equation}
f'(Y^*) = 1 - \left[\sin{\alpha} + 2\sqrt{(Y^* - \cos\alpha)(E_s - Y^*)}\right]\frac{d}{d\:Y_i}\Delta\theta_i\Big\rvert_{Y^*}
\label{stability3n}
 \end{equation}
by using $\Delta \theta_{i} = 2\alpha$ for symmetric contact phases. Later in the paper, we take the notation  $\frac{d}{d\:Y_i}\Delta\theta_i\Big\rvert_{Y^*}=d_i\Delta \theta^*$.
Since the expression in square brackets in  (\ref{stability3n}) always remains positive,  then (\ref{stability2n})
transforms into the following condition for the angle swept during stance
\begin{equation}
d_i\:\Delta \theta^* \in \left(0,\frac{2}{\sin{\alpha} + 2\sqrt{(Y^* - \cos{\alpha} )(E_s - Y^*)}}\right),
\label{stability33n}
 \end{equation}
 which shows that
higher or lower apex heights must be compensated appropriately
 larger or smaller amount of angular sweep $\Delta\theta_i$. 
 However, for stability reasons, the positive derivative $d_i\:\Delta \theta^*$ cannot be greater than
  $2\left[\sin{\alpha} + 2\sqrt{(Y^* - \cos{\alpha})(E_s - Y^*)}\right]^{-1}$. 
Then we substitute 
the equations fixed point $Y^*$  (\ref{apexheight3}) and (\ref{apexheight3nowy})
 satisfying the condition (\ref{energyminnowy}) and denoting $2\sqrt{(Y^* - \cos{\alpha})(E_s - Y^*)}$ by $D\left(\alpha,E_s,\theta_d^*\right)$ we get
 \begin{equation}
D\left(\alpha,E_s,\theta_d^*\right)=\theta^*_d\cos{2\alpha} \sqrt{2E_s-(\theta^*_d)^2-2\cos{\alpha}}+\sin{2\alpha}\left[E_s-(\theta^*_d)^2-\cos{\alpha}\right].
\label{D_a,Es,thd}
 \end{equation}

Next, from (\ref{angular333}) it follows that
\begin{equation}
\begin{split}
	& d_i\Delta \theta^*= \frac{d_i\theta_d^*}{\theta^*_d}\left[
2\alpha +\frac{\pi^2\sin{\alpha}}{2K}\right] \\
& + \frac{2\theta^*_d}{K} \left[2 d_i L_d^* -  2\frac{\theta^*_d}{L^*_d} d_i \theta_d^* +\frac{(\theta_d^*)^2-\cos{\alpha}}{(L^*_d)^2}d_i L_d^*\right],
\end{split}
\label{}
 \end{equation}
 where $K$ is given by (\ref{approx444}).
Resolving $\partial_i \theta_d^*$ (see (\ref{ThetaLdUV2}))
 \begin{equation}
 d_i \theta_d^*=\frac{d}{d\: Y_i}\theta_d\Big\rvert_{Y^*}= -\frac{1}{\sqrt{2}}\left(\frac{\cos \alpha}{\sqrt{E_s-Y^*}}+\frac{\sin \alpha}{\sqrt{Y^* - \cos{\alpha}}}\right), 
\label{stability444n}
 \end{equation}
  and applying (\ref{map3333}) we have
\begin{equation}
 d_iL_d^*=\frac{d}{d\: Y_i}L^*_d\Big\rvert_{Y^*}= -\frac{\theta^*_d}{L^*_d}\, d_i\theta_d^*, 
\label{stability555n}
 \end{equation}
 and from (\ref{stability555n}) can be further deduced 
 \begin{equation}
\begin{split}
	& d_i\Delta\theta^*=-2d_i \theta_d^* \times \\
	&  \times \left[\frac{(\theta_d^*)^2}{K\sqrt{2E_s-(\theta^*_d)^2-2\cos{\alpha}}}\left(4+\frac{(\theta^*_d)^2-\cos{\alpha}}{2E_s-(\theta^*_d)^2-2\cos{\alpha}}\right) -\frac{1}{\theta^*_d}\left(
\alpha +\frac{\pi^2\sin{\alpha}}{4K}\right) \right]=\\
	& =\frac{2}{D\left(\alpha,E_s,\theta_d^*\right)} \times \\
&  \times \left[\frac{(\theta_d^*)^2}{K}\left(4+\frac{(\theta^*_d)^2-\cos{\alpha}}{2E_s-(\theta^*_d)^2-2\cos{\alpha}}\right) -\frac{\sqrt{2E_s-(\theta^*_d)^2-2\cos{\alpha}}}{\theta^*_d}\left(
\alpha +\frac{\pi^2\sin{\alpha}}{4K}\right) \right].
\end{split}
\label{stability4nowy}
 \end{equation}
with $D\left(\alpha,E_s,\theta_d^*\right)$ given by (\ref{D_a,Es,thd}). 

\begin{thm} \textnormal{(On the stability of symmetric solutions).}
 Let us assume that Theorem \ref{th1} holds. If additionally, the set of parameters $\{\alpha, E_s, \theta^*_{d}\}$ satisfies the condition 
\begin{equation}
d_i\Delta\theta^* \in \left(0,\frac{2}{\sin{\alpha} + D\left(\alpha,E_s,\theta_d^*\right)}\right)
\label{stability3333n}
 \end{equation}
with $d_i\Delta\theta^*$ given by the formula (\ref{stability4nowy}),
 then $Y^*$ given by (\ref{apexheight3}) or (\ref{apexheight3nowy}) is a stable fixed point.
 \label{th2}
 \end{thm} 
 \begin{rem}
When the angular velocity at touch-down for the fixed point $Y^*$ equals $\theta_d^*=\sqrt{\cos\alpha}$, then the condition (\ref{stability3333n}) from Theorem \ref{th2} takes the form 
\begin{equation}
d_i\Delta\theta^* \in \left(0,\frac{2}{\sin{\alpha} + D(\alpha,E_s)}\right),
\label{stability4444n}
 \end{equation}
 where from (\ref{D_a,Es,thd})
 \begin{equation}
 D(\alpha,E_s)=\cos{2\alpha} \sqrt{\cos{\alpha}(2E_s-3\cos{\alpha})}+\sin{2\alpha}\left[E_s-2\cos{\alpha}\right].
\label{D(C=0)}
 \end{equation}
 If we additionally assume that $K$ is large enough to be replaced by $\widetilde{K}=\frac{\pi^2 \cos{\alpha}}{4\alpha^2}$ (see (\ref{napprox555})), then the formula for $d_i\Delta\theta^*$ given by (\ref{stability4nowy}) takes the form
  \begin{equation}
d_i\Delta\theta^* =\frac{2\alpha}{D(\alpha,E_s)} \times \left[\frac{16\alpha}{\pi^2} -\sqrt{
\frac{2E_s-3\cos{\alpha}}{\cos{\alpha}}}\left(
1 +\alpha\tan{\alpha}\right) \right].
\label{stability5n}
 \end{equation}
 \label{rem2}
 \end{rem}

\subsection{ Stability regions for special case $\theta_d^*=\sqrt{\cos\alpha}$ }
\label{a0}
On the one hand, the special case is a mathematical simplification compared to the general one. On the other hand, it reflects a typical running speed in humans.  For small angles of
attack $\alpha$, the horizontal Froude number $X'$ relates to the angular
velocity at
touch-down with $X' \approx \theta_d$. Hence the case $C(\epsilon) = 0$  describes running  for which the horizontal velocity $X'$ is approximately 1. For a leg length of one meter,  this means human jogging speed about $3.13\;m/s$, 
 which exactly corresponds to the initial velocity for a minimum energy of 1.5 and an initial apex height of 1.

 Two curves: lower, when $f'(Y^*)=-1\;$ i.e. $\;d_i\Delta\theta^* =2[\sin{\alpha} + D(\alpha,E_s)]^{-1}$,  and upper, when $f'(Y^*)=1\;$ i.e. $\;d_i\Delta\theta^* =0$ (see \ref{stability2n}) are the region constraints:  $E_s(\alpha)$ and $K(\alpha,E_s)$,  consisting of a combination of parameters leading to a stable fixed point. Due to Remark \ref{rem2} we obtain the following $E_s(\alpha)$ - region:
 \begin{equation}
E_s^-<E_s^{min} \quad <\quad E_s \quad < \quad E_s^+=\frac{3}{2}\cos{\alpha}+\frac{128\alpha^2\cos{\alpha}}{\pi^4(1+\alpha\tan{\alpha})^2},
\label{energy222new}
 \end{equation}
 where $E_s^-$ and $E_s^+$ are determined from (\ref{stability4444n}) and (\ref{stability5n}), and $E_s^{min}=\frac{1}{2\cos{\alpha}}+\cos{\alpha}$ is given by the condition (\ref{energy222}), and also fixed $\alpha \in \left(\pi/36,\pi/6\right)$. It is easy to verify that $E_s^-<E_s^{min}$ for each parameter $\alpha$ within the range under consideration. Thus, the lower energy limit of the system is $E_s^{min}$. Moreover, using (\ref{approx444}) and (\ref{delta3}) we can find the following $K(\alpha, E_s)$ - region  for stable solutions
 \begin{equation}
K^{min} \quad < \quad K  \quad< \quad K^+,
\label{stiffnessnew}
 \end{equation}
where $K^{min}$ is given by (\ref{stiffness222}) and (\ref{deltamin222}), $K^+=K(\alpha,E_s^+)$, and also fixed $\alpha \in \left(\pi/36,\pi/6\right)$.

\begin{figure}
		\centering 
		\includegraphics[width=0.49\textwidth]{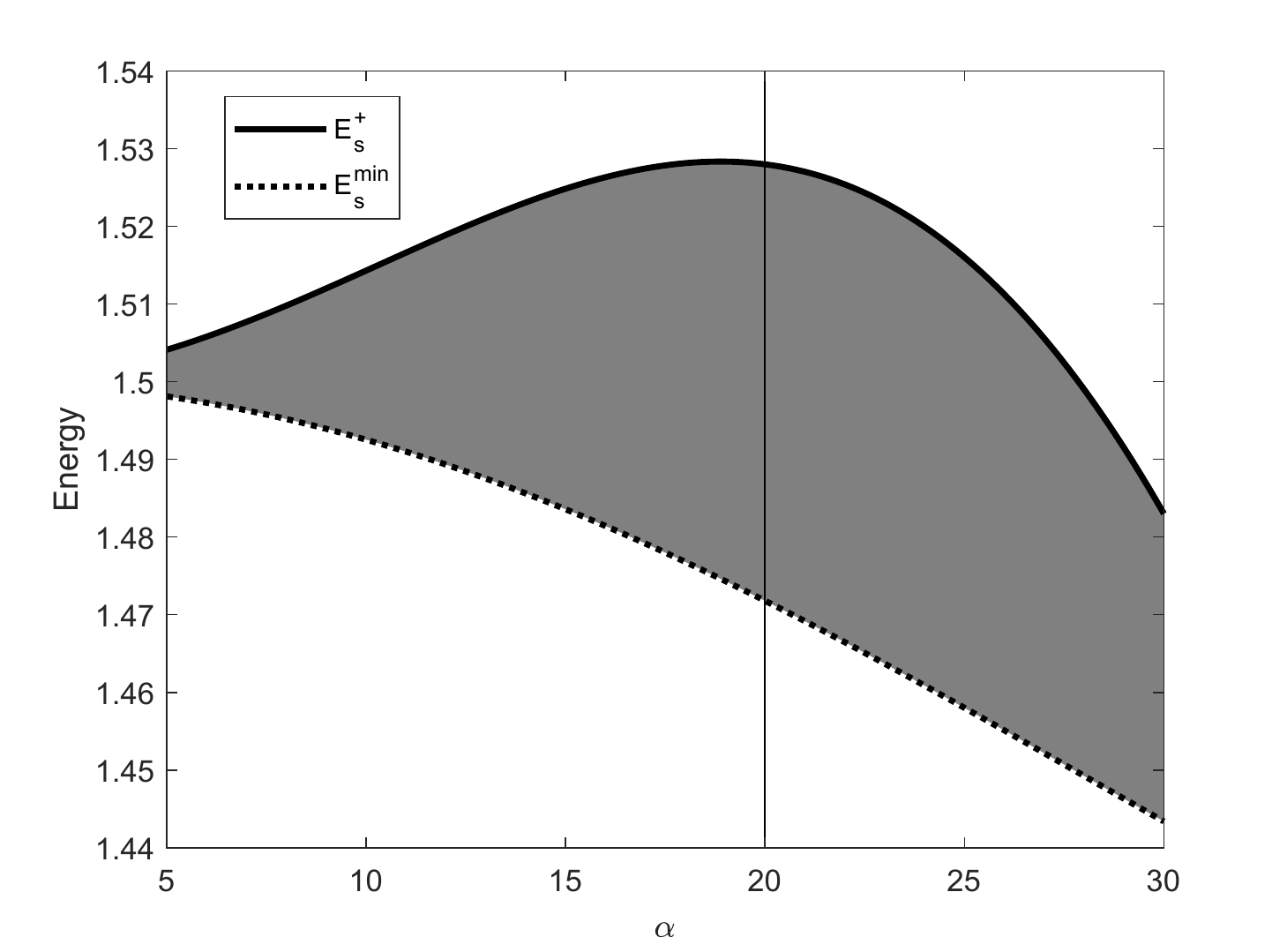}
		\includegraphics[width=0.49\textwidth]{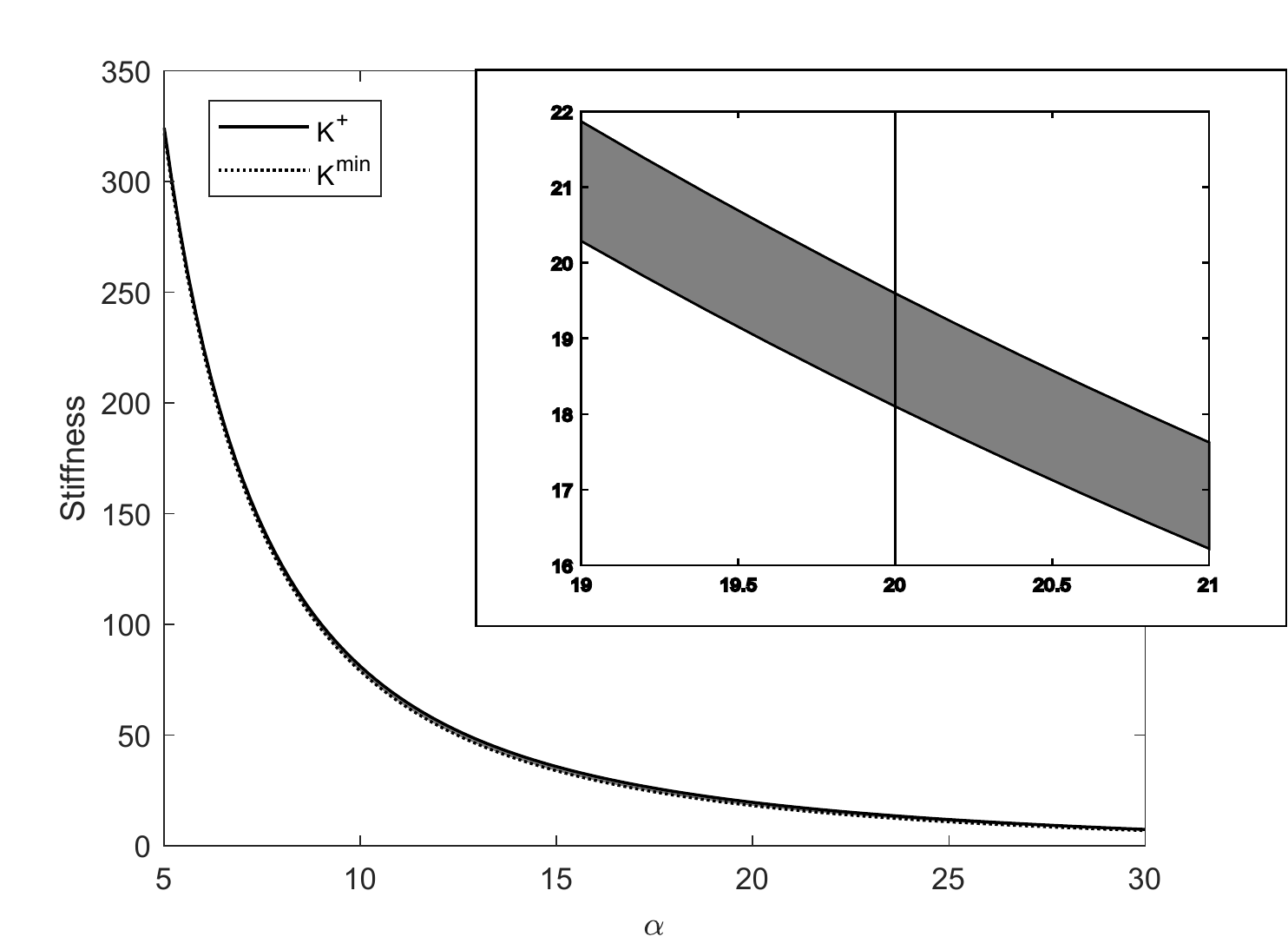}
		\caption{\small Parameter interdependence for fixed point solutions with $\theta_d^*=\sqrt{\cos\alpha}$. The stability regions: $E_s(\alpha)$ (on the left) and $K(\alpha,E_s)$ (on the right) for spring-mass running predicted by the analytical approximations (\ref{energy222new}) and (\ref{stiffnessnew}) for varying  angles of attack $\alpha$ from  $5^{\circ}(=\pi/36)$ to $30^{\circ}(=\pi/6)$.}
		\label{Fig4}
\end{figure}

Based on these results (see (\ref{energy222new}) and (\ref{stiffnessnew})), parameter combinations leading to stable fixed points are indicated in Fig. \ref{Fig4} as dark areas. For example, an angle of
attack $\alpha=\pi/9$, marked in Fig. \ref{Fig4} with a vertical line,   necessitates a
minimum system energy $E_s^{min}=1.4718$ (shown in Fig. \ref{Fig4}) and the lower
stability constraint corresponds to a system energy
$E_s^{-}=1.411$ below this minimum (not shown in Fig. \ref{Fig4}). On the other hand,  the
system energy related to the upper constraint
$E_{s}^{+}=1.528$ (shown in Fig. \ref{Fig4}) still exceeds the critical level. Moreover, on the right side of Figure \ref{Fig4}, we can check that for $\alpha=\pi/9$: $$18.1013=K^{min} \quad < \quad K \quad < \quad K^+=19.5960.$$

Taking into account the assumption of small values of the angle $\alpha$, the range of the approximate solution is always bound to a
spring stiffness exceeding physiologically reasonable
values. Since  angles of attack  less than $\pi/18$ are unlikely to happen in a real human run with jogging speed, then the stiffness $K$ determined from the symmetrical model (see (\ref{approx444}) and (\ref{delta3})) cannot be a reliable physiological value. For  example, taking $\alpha=\pi/36$  and
$E_s=1.5$, the dimensionless stiffness $K(\pi/36,1.5)$ is 322 (see Fig. \ref{Fig4}). Compared to experiments, where typical stiffness values are in the range of 40 - 80 (e.g. \cite{wrobel}), the obtained $K = 322$ is from 4 to 8 times greater.

\begin{figure}
		\centering 
		\includegraphics[width=0.9\textwidth]{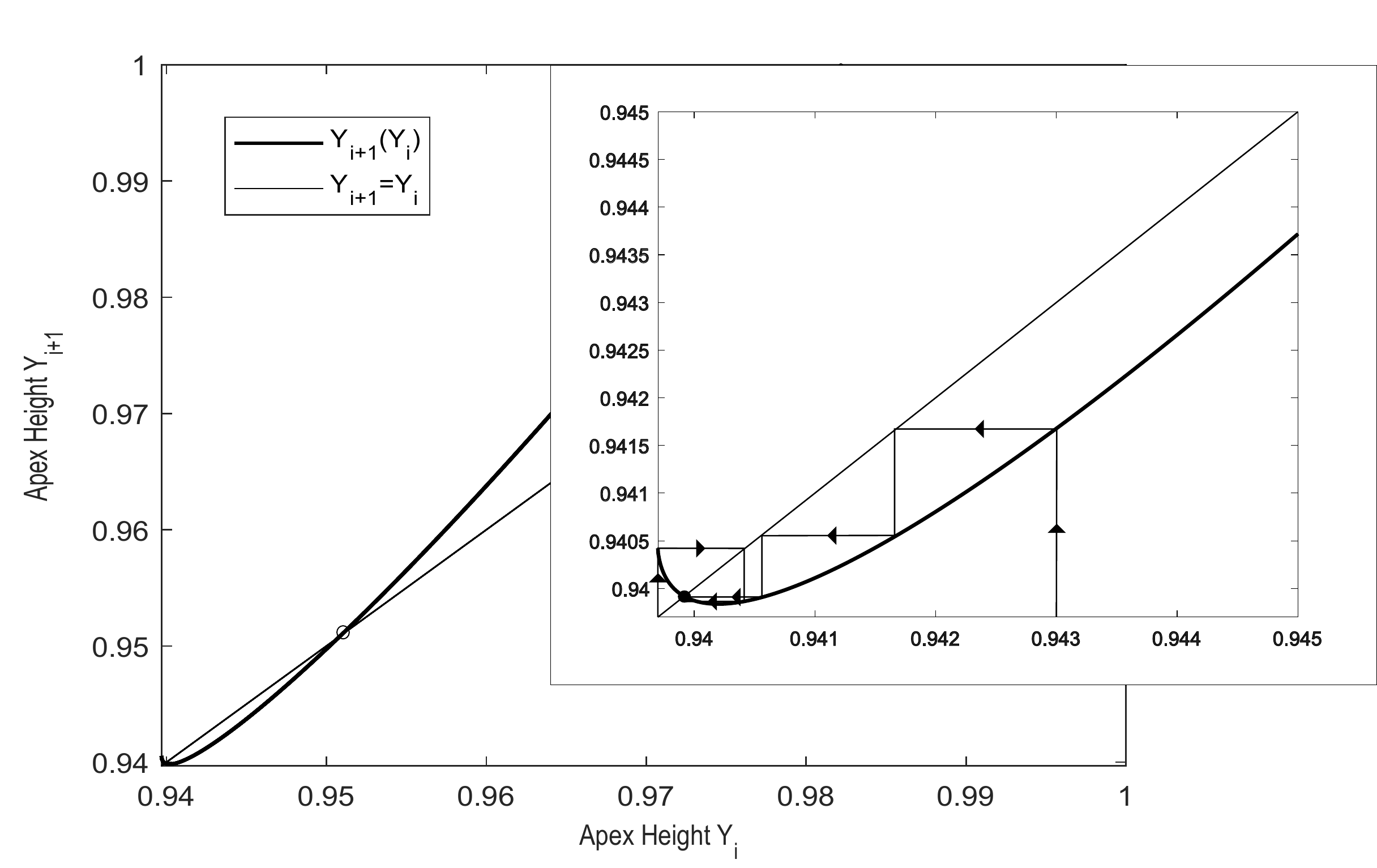}
		\caption{\small Stability of spring-mass running. The return map function $Y_{i+1}(Y_i)$ is given by (\ref{map3}) for the parameter set:  $\alpha=\pi/9$, $E_s=1.48$ and $K=18.3575$ with marked the stable
fixed point (black circle) and the unstable
fixed point (white circle).}
		\label{Fig5}
\end{figure}

The return map function $Y_{i+1}(Y_i)$ (see (\ref{map3})) in Fig. \ref{Fig5}  is shown for the parameter set  $\alpha=\pi/9$, $E_s=1.48$ and $K=K(\pi/9,1.48)=18.3575$, which belongs to the calculated regions of parameter combinations producing stable fixed point solutions. As predicted, the return
map has a stable fixed point $Y^*$ (see (\ref{apexheight3}) or (\ref{apexheight3nowy})) attracting neighboring apex
states within a few steps (arrow traces in the magnified region in Fig. \ref{Fig5}).
The stable fixed point  $Y^*$ is calculated from the formula (\ref{apexheight3})  or (\ref{apexheight3nowy}) and for the parameter set: $\alpha=\pi/9$, $E_s=1.48$, $\theta_d^*=\sqrt{\cos{(\pi/9)}}$, and  $L_d^*=\sqrt{2\times 1.48-3\times\cos{(\pi/9)}}$, it is $0.9399$, exactly as shown in Fig. \ref{Fig5}. Furthermore, the return map is characterized by an additional,
unstable fixed point representing the upper limit of the basin of
attraction of the stable one. 
Here, the basin of attraction contains all apex
heights from the landing height $Y_i=\cos{\alpha}$ to the
second, unstable fixed point (white circle).  Both fixed points
are merely observations from plotting equation (\ref{map3}). However, from Fig. \ref{Fig5} it is obtained that, if the system energy
is not adequately selected ($E_s>E_s^+$), the fixed point given by (\ref{apexheight3})  or (\ref{apexheight3nowy}) is unstable.

\subsection{Stability regions for generalized symmetrical cases}
\label{general}
\noindent
In this section, we will deviate from the case of $\theta_d^*=\sqrt{\cos\alpha}$  and present the regions of stability in a general approach. The upper and lower limits of the stability regions, denoted in the figures by a solid and a  dashed line, respectively, are obtained  from Theorem \ref{th2}.

Fig. 6 shows the stability regions of angular velocity $\theta_d^*$ depending on the changing energy (with a given $\alpha$) and the changing angle of attack (with a given $E_s$).
For a fixed angle of attack, the angular velocity $\theta_d^*$ increases with the increase in energy, while for a fixed energy, the angular velocity  $\theta_d^*$ decreases with the increase in the angle of attack.

\begin{figure}
		\centering 
		\includegraphics[width=0.49\textwidth]{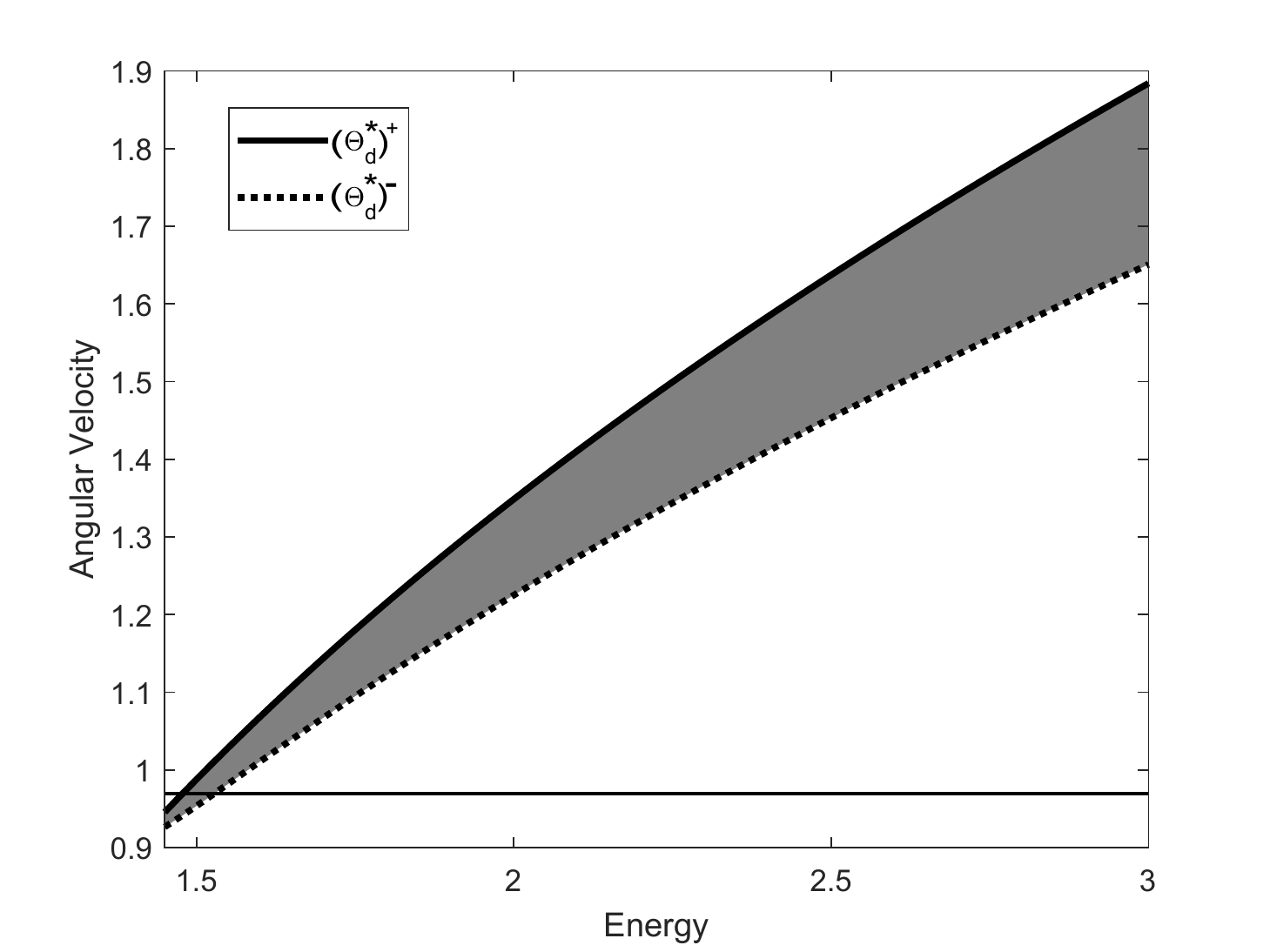}
		\includegraphics[width=0.49\textwidth]{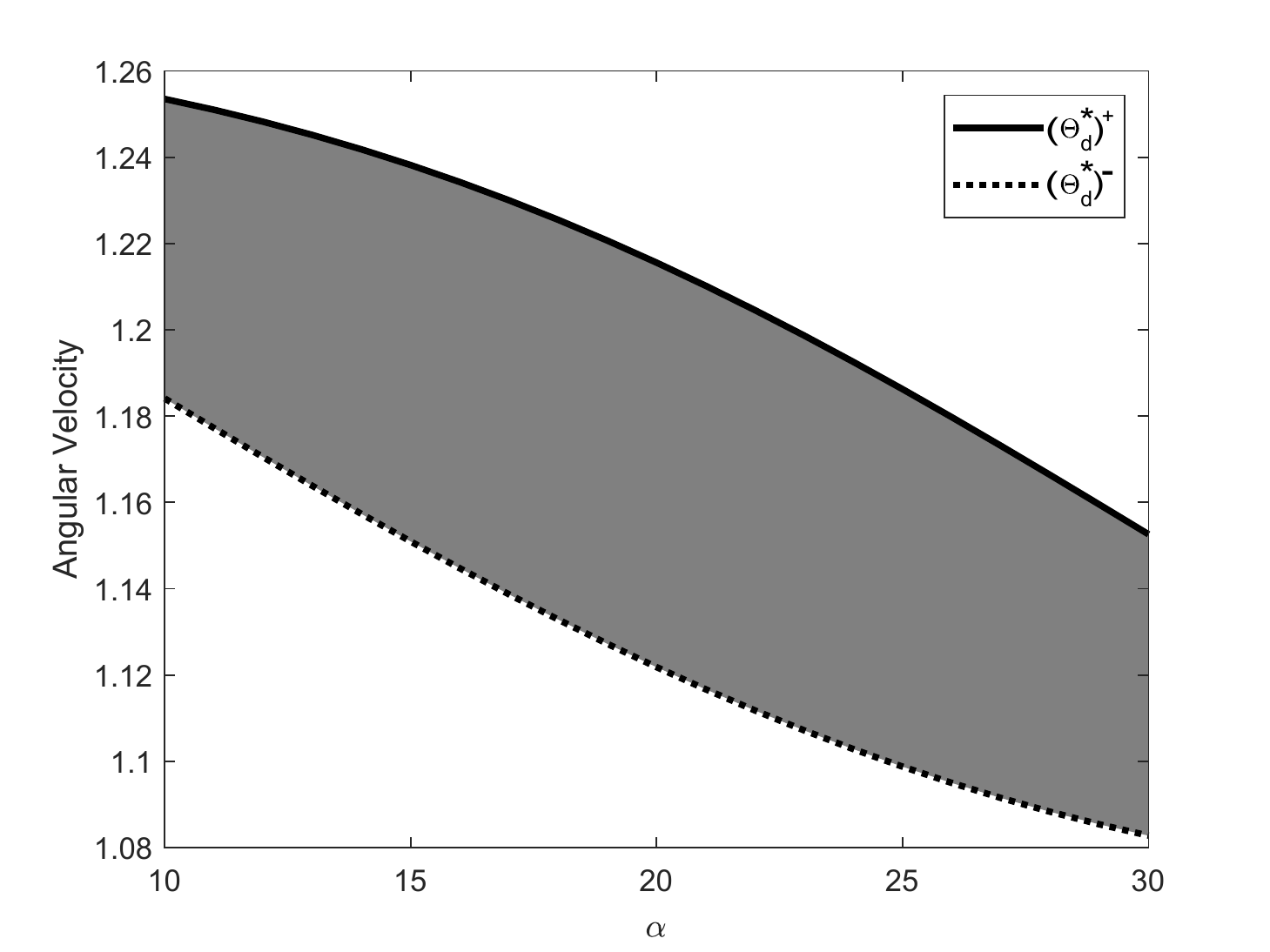}
		\caption{\small On the left, the stability region of  the angular velocity $\theta_d^*(E_s)$ with $\alpha=\pi/9$ for varying energy $E_s$ from $1.45$ to $3$. The horizontal line shows the case of $\theta_d^*=\sqrt{\cos\alpha}=0.9694$.
		On the right, the stability region of  the angular velocity $\theta_d^*(\alpha)$ with $E_s=1.8$ for varying  angles of attack $\alpha$ from $10^{\circ}(=\pi/18)$ to $30^{\circ}(=\pi/6)$.
		The upper and lower limits are $\left(\theta_d^*\right)^+$ and $\left(\theta_d^*\right)^-$, respectively.}
		\label{Fig6}
\end{figure}
For sprints, running technique changes significantly. This is due to higher $\theta_d^*$ values and smaller $\alpha$ values. While the angles of attack are smaller than in jogging, the take-off angles are significantly larger. This introduces an asymmetry for which the model is not suitable, so it does not make sense to analyse the stability for this running technique. In this section physiologically valid values of leg stiffness are considered (i.e. $\alpha$ greater than $\pi/18$).

\begin{figure}
		\centering 
		\includegraphics[width=0.49\textwidth]{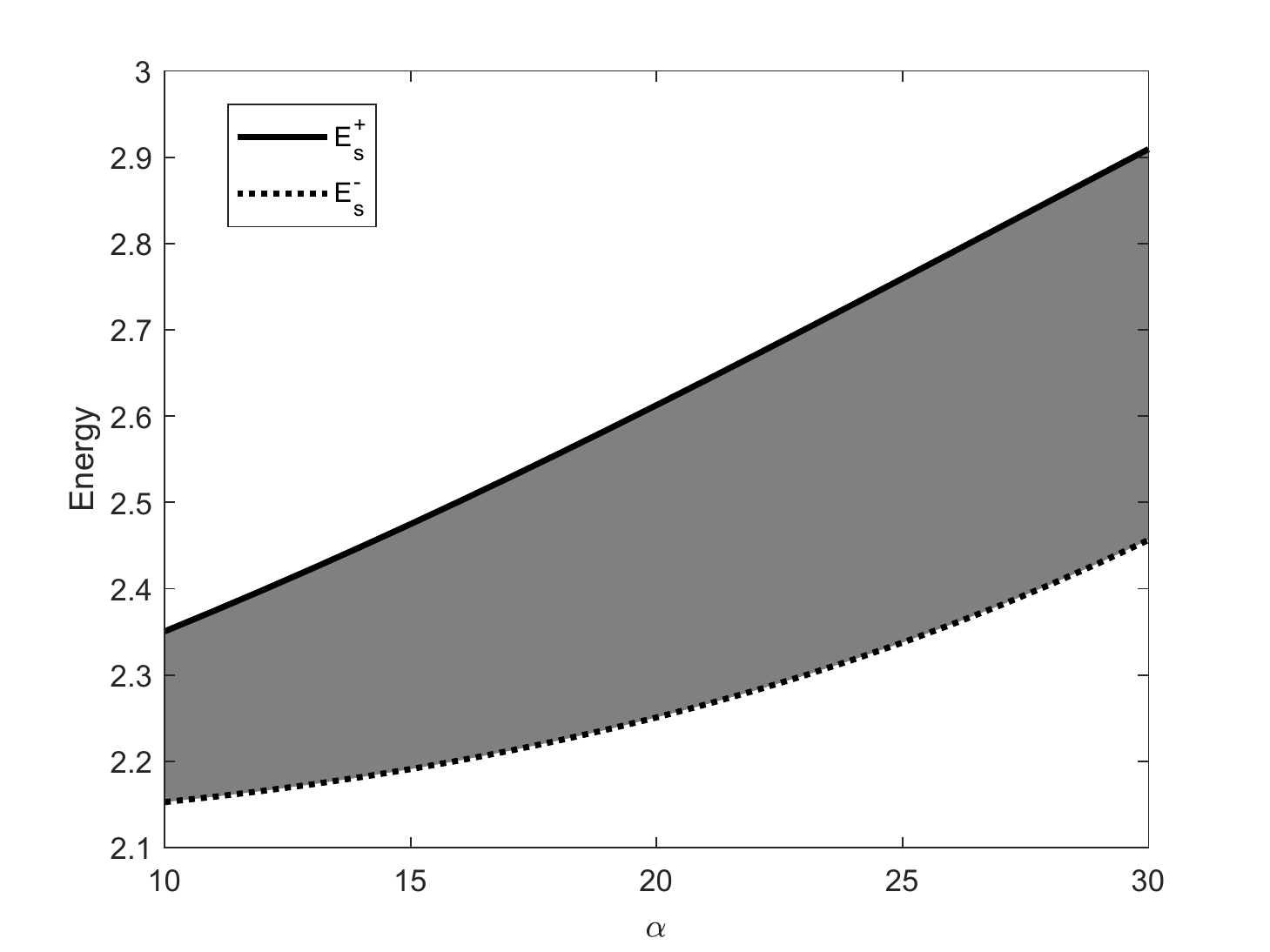}
		\includegraphics[width=0.49\textwidth]{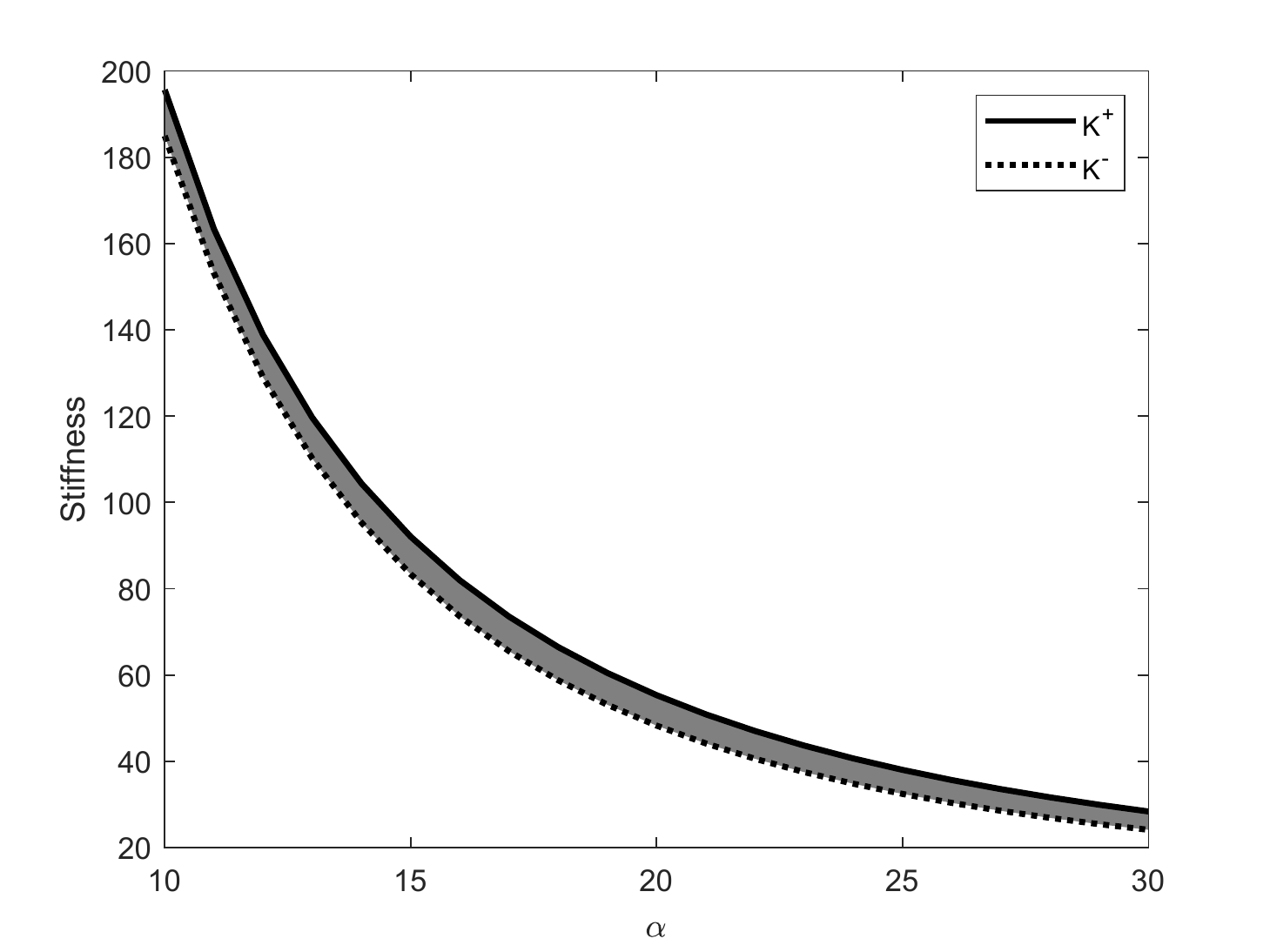}
		\caption{\small 
		On the left, the stability region of  the energy $E_s(\alpha)$, while on the right,  the stability region of the dimensionless stiffness $K(\alpha, E_s)$, determined by (\ref{approx444}) and (\ref{delta3}),
		with  $\theta_d^*=1.5$ for varying  angles of attack $\alpha$ from $10^{\circ}(=\pi/18)$ to $30^{\circ}(=\pi/6)$.
		The upper and lower limits are respectively $E_s^+$ and $E_s^-$ on the left, and $K^+=K(\alpha, E_s^+)$ and $K^-=K(\alpha, E_s^-)$ on the right. }
		\label{Fig7}
\end{figure}

Fig. \ref{Fig7} illustrates the stability regions of energy $E_s$ and stiffness $K(\alpha,E_s)$, defined by (\ref{approx444}) and (\ref{delta3}), for the changing angle of attack in the case of a fixed angular velocity.
As can be seen in Fig. \ref{Fig7}, the energy $E_s$ increases with increasing angle of attack at constant angular velocity, while the stiffness $K$ of course decreases.

\begin{table}[t]
\small
\centering
\begin{tabular}{lccccc|ccccc}
\toprule
  &\multicolumn{5}{c|}{$\alpha \quad\quad (\theta_d^* =1.18)$} & \multicolumn{5}{c}{$\theta_d^* \quad\quad (\alpha=20\pi/180)$}  \\

 & $\frac{10\pi}{180}$ & $\frac{15\pi}{180}$ & $\frac{20\pi}{180}$ & $\frac{25\pi}{180}$ & $\frac{30\pi}{180}$ & $0.95$ & $1$ & $1.05$ & $1.10$ & $1.15$  \\
  \midrule 
  $E_s^{min}$ & 1.703 & 1.712 & 1.728 & 1.754 & 1.794 & 1.451 & 1.506 & 1.564 & 1.625 & 1.688  \\
  \hline
  $K^{min}$ & 85.06 & 37.94 & 21.60 & 14.21 & 10.37 & 17.79 & 18.60 & 19.42 & 20.26 & 21.10  \\
  \hline
  $\Delta L_{MAX}$ & 0.0449 & 0.0807 & 0.1216 & 0.1655 & 0.2105 & 0.1415 &0.1367 & 0.1321 & 0.1279 & 0.1239   \\
\bottomrule
\end{tabular}
\caption{Interdependence of parameters $E_s^{min}$, $K^{min}$ and  $\Delta L_{MAX}$, given by formulas   (\ref{energymin}),  (\ref{stiffnessaprox}), (\ref{delta5}) and  (\ref{maksspring}) in the spring-mass running.}
\label{tab}
\end{table}

If the angular velocity $\theta_d^*$ is constant, the energy increases with increasing angle of attack $\alpha$. During the contact phase, for a larger angle of attack $\alpha$, more energy is required for the mass point to stabilize at the same level. On the other hand, for a fixed $\alpha$, more energy is needed to increase the angular velocity $\theta_d^*$. The stiffness of the spring depends both on the angle $\alpha$ and the energy  $E_s$. The greater the energy, the greater the stiffness, while the stiffness decreases with increasing angle of attack. Let us note, that the effect of $\alpha$ on leg stiffness is much stronger and a typical situation is a decrease in stiffness with an increase in both the angle of attack and energy. It is easy to conclude that with greater leg stiffness (small $\alpha$), the spring deflects less. In Section 3, we observed  that small deflection is also influenced significantly by large angular velocity $\theta_d^*$ (parameter $A(\epsilon)$ is positive). 

The stability of the model enforces the symmetry of solutions, which for the small $\alpha$, required for approximate solutions, gives unreasonable values of the running parameters. Thus,
the compromise requires to analyze the stability for the attack angles in the range from $\pi/18$ to $\pi/6$.

\subsection{Transcritical bifurcation}
\label{Bifurcation}

\begin{figure}
		\centering 
		\includegraphics[width=0.7\textwidth]{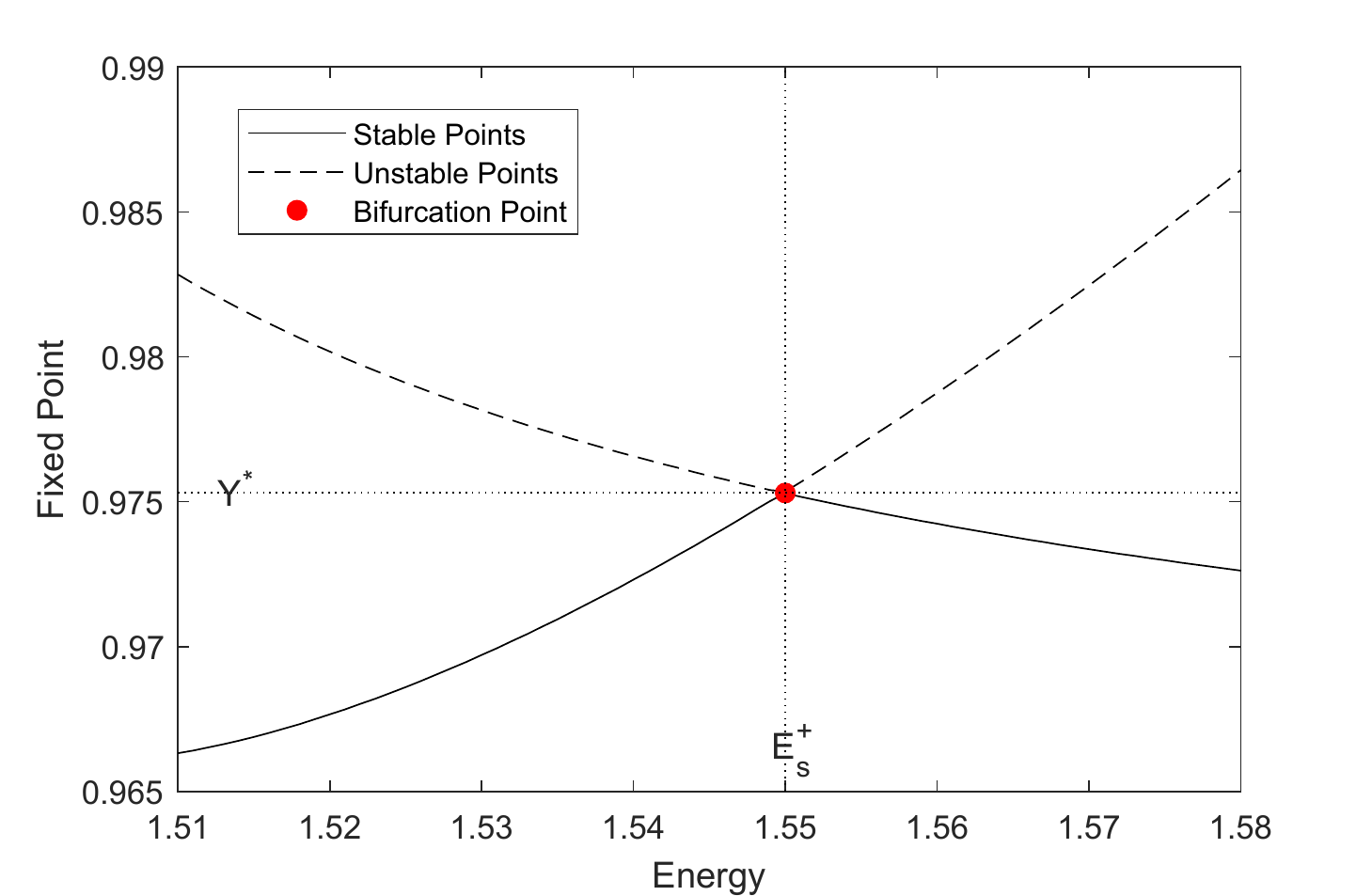}
		\caption{\small One-parameter bifurcation diagram of (\ref{map3}). Note transcritical bifurcation of fixed point $Y^* = 0.9754$ at $E_s^+=1.550$. The angle of attack $\alpha$ is set to $\alpha=\pi/12$, while the stiffness $K$ obtained from (\ref{approx444}) increases from $34.77$ for $E_s = 1.51$ to $38.73$ for $E_s = 1.58$.}
		\label{Fig8}
\end{figure}

\noindent In the current section, we will present numerical evidence for the occurrence of a transcritical bifurcation in our mapping given by equation (\ref{map3}). In particular, by considering $E_s$ as a bifurcation parameter, we will show that fixed point solutions collide in a transcritical bifurcation and exchange stability. Let us focus in Fig.~\ref{Fig8}, which is a one-parameter bifurcation diagram, where we depict the existence and stability of fixed point solutions versus energy. The solid line  in the figure indicates stable solutions, and the dashed line unstable ones.        

Let us start from the left of the figure with the lowest energy level $E_s = 1.51$. At this value of the bifurcation parameter, there is a pair of fixed point solutions, the one on the lower branch is stable and the one on the upper branch is unstable. Increasing the value of energy parameter $E_s$, these two solutions move closer together until they collide for $E_s^+ = 1.550$ with $Y^* = 0.9754$. Increasing $E_s$ further, past the value of $E_s = E_s^+$ we can see that the stable branch continues past the bifurcation point as an unstable branch, and the unstable one as the stable branch. We thus have a typical scenario of a transcritical bifurcation. 
\begin{figure}
		\centering 
		\includegraphics[width=0.7\textwidth]{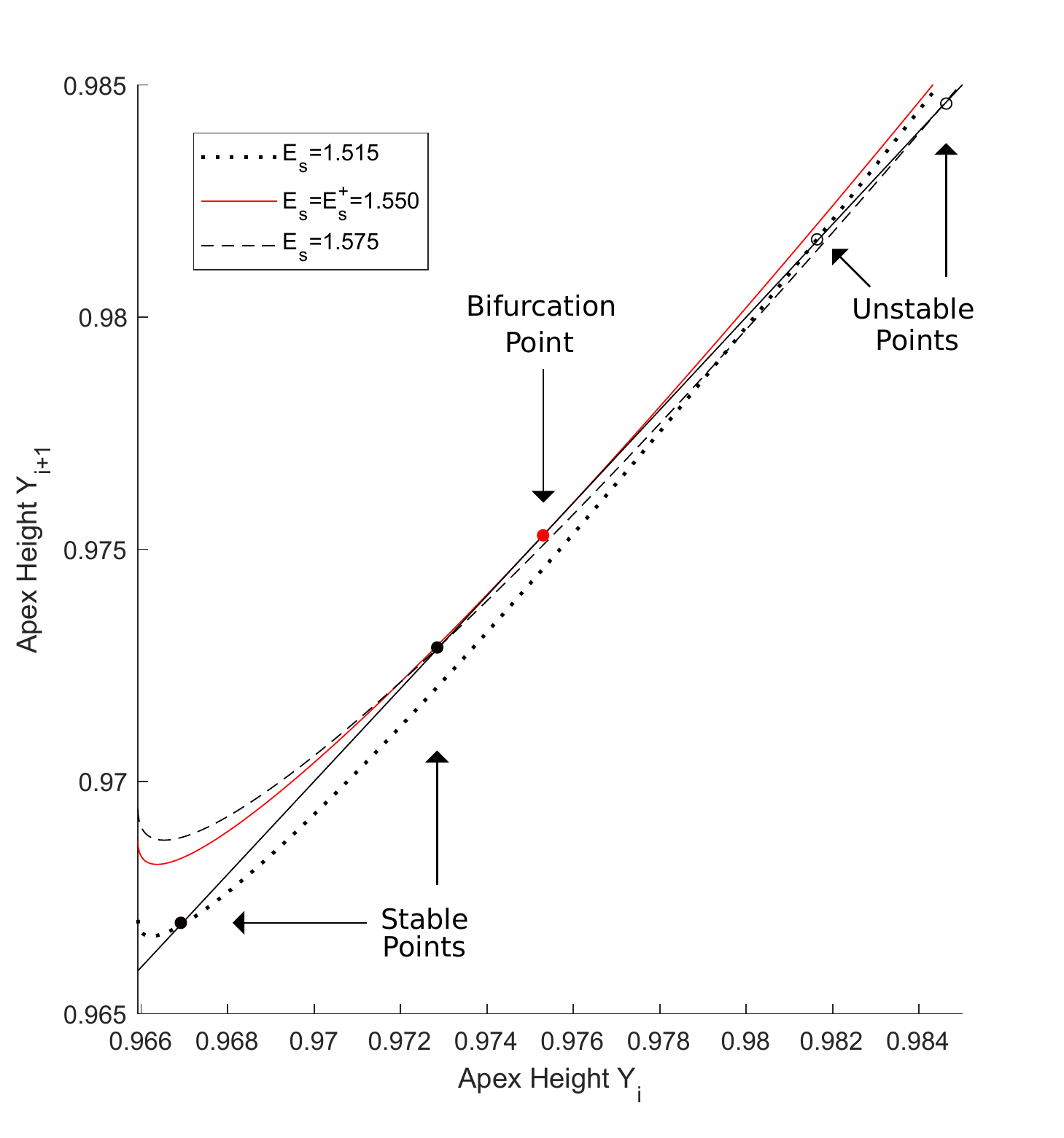}
		\caption{\small Numerical depiction of transcritical bifurcation occurring under variation of bifurcation parameter $E_s$ in mapping (\ref{map3}). The angle of attack $\alpha$ is set to $\alpha=\pi/12$, while the stiffness $K$ obtained from (\ref{approx444}) is $35.07$ for $E_s = 1.515$,  $36.75$ for $E_s=E_s^+ = 1.550$ and $37.94$ for $E_s=1.575$.}
		\label{Fig9}
\end{figure}
To see this further, in Fig.~\ref{Fig9}, we numerically depict mapping (\ref{map3}) for three distinct values of the bifurcation parameter. Namely for $E_s < E_s^+$, then for $E_s = E_s^+$ and, finally for $E_s > E_s^+$. Let us consider first the graph shown by black dotted curve for $E_s = 1.515 < E_s^+$. We can see in the graph two fixed points with the stable (black solid dot) one below the unstable one (black empty circle), which agrees with the bifurcation diagram presented earlier in Fig.~\ref{Fig8}. Increasing $E_s$ to $E_s = E_s^+ = 1.550$ the two fixed points collide - see the red solid curve in the figure. Notice that the graph is quadratically tangent to line $Y_{i + 1} = Y_i$ at the bifurcation point $Y^*$ (red solid dot). Finally, increasing $E_s$ to $E_s = 1.575 > E_s^+$ allows us to see the exchange of stability of the fixed points - see the black dashed curve in Fig.~\ref{Fig9} - the unstable fixed point (black empty circle), say $Y_u$ still lies  above fixed point $Y^*$, that is $Y_u > Y^*$, and the stable one (black solid dot), say $Y_s$, lies below $Y^*$, that is $Y_s < Y^*$. From these three graphs we can see that the defining and nondegeneracy conditions for a transcritical bifurcation are satisfied. Namely, we observe a quadratic tangency of the mapping at fixed point $Y^*$ at $E_s = E_s^+$, and unfolding with respect to parameter $E_s$, which shows the exchange of stability of the fixed points.     

\section{Conclusions}
\label{conclusion}
\noindent
The analysis presented in the paper shows that we can construct a useful and relatively simple model as an approximation to the more complex spring-mass description of running. In \cite{geyer}, the authors used a similar strategy to reduce the spring-mass model to a one-dimensional map. However, their method of constructing a one-dimensional mapping as a reduced approximate model is different than ours. That is, we are led by rigorous application of scaling and perturbation analysis conducted in \cite{plocin}, which allows us to neglect several small terms and remain asymptotically consistent. This rigorous approximation leads to an integrable system that can be thoroughly analysed. In particular, we have proved that under some natural conditions the apex to apex mapping has a stable fixed point which concurs with realistic situations. We conduct stability analysis of this fixed point to determine the regions of its existence in parameter space. We also determine the condition for the existence of fixed points expressed as an inequality linking the system's energy with the angle of attack and angular velocity. This result is stated in Theorem 1. We also perform a numerical continuation of fixed points with respect to energy. This has allowed us to discover that the stable fixed point undergoes a transcritical bifurcation. A thorough investigation of this behaviour as well as the analysis of asymmetric solutions, which are of special interest in the modelling of running, is the subject of our further studies.  

\printbibliography

\end{document}